\documentclass[12pt]{amsart}
\usepackage{geometry}
\usepackage{amssymb}
\usepackage{amsmath}
\usepackage{amsthm}
\usepackage{graphicx}
\usepackage{bm}
\usepackage[usenames,dvipsnames]{xcolor}
\usepackage{tikz}
\usepackage{subcaption}
\usepackage{float}
\usepackage{color}
\usepackage{hyperref}

\usepackage{appendix}
\usepackage{array}
\usepackage{footnote}
\usepackage{booktabs}
\usepackage{hhline}
\makesavenoteenv{tabular}

\usepackage{algpseudocode}
\usepackage{algorithm}

\numberwithin{equation}{section}

\makeatletter
\g@addto@macro{\endabstract}{\@setabstract}
\newcommand{\authorfootnotes}{\renewcommand\thefootnote{\@fnsymbol\c@footnote}}%
\makeatother

\newcommand{\R}{\mathbb{R}}
\DeclareMathOperator*{\argmin}{arg\,min}

\begin{document}

 \begin{center}
 \large
   \textbf{Backtracking strategies for accelerated descent methods with smooth composite objectives} \par \bigskip \bigskip
   \normalsize
  \textsc{Luca Calatroni}\textsuperscript{$\dagger$}, \textsc{Antonin Chambolle} \textsuperscript{$\dagger$}\let\thefootnote\relax
\footnote{
\textsuperscript{$\dagger$} Centre de Math\'ematiques Appliqu\'ees (CMAP), \'Ecole Polytechnique CNRS, 91128, Palaiseau Cedex, France.
}
\let\thefootnote\relax\footnote{
\hspace{3.2pt} e-mail: \href{mailto:luca.calatroni@polytechnique.edu}{\nolinkurl{luca.calatroni@polytechnique.edu}},
 \href{mailto: antonin.chambolle@cmap.polytechnique.fr}{\nolinkurl{antonin.chambolle@cmap.polytechnique.fr}}}

\end{center}

\begin{abstract}
We present and analyse a backtracking strategy for a general Fast Iterative Shrinkage/Thresholding Algorithm which has been proposed in \cite{ChambollePock2016} for strongly convex composite objective functions. Differently from classical Armijo-type line searching, our backtracking rule allows for local increasing and decreasing of the descent step size (i.e. proximal parameter) along the iterations. We prove accelerated convergence rates and show numerical results for some exemplar imaging problems.
\end{abstract}



\section{Introduction}  \label{sec:intro}

The concept of \emph{acceleration} of first-order optimisation methods dates back to the seminal work of Nesterov \cite{Nesterov1983}. For a proper, convex, l.s.c. function $F:\mathcal{X}\to\R\cup\left\{\infty\right\}$ defined on a Hilbert space $\mathcal{X}$ with Lipschitz gradient with constant $L>0$, solving the abstract optimisation problem 
\begin{equation} \label{opt:problem}
\min_{x\in\mathcal{X}}~F(x)
\end{equation}
by means of an accelerated iterative method means improving the convergence rate $O(1/k)$ achieved after $k\geq 1$ iterations of standard gradient descent methods in order to (almost) match the universal lower bound of $O(1/k^2)$ holding for any function such as $F$. In the smoother case, i.e. when $F$ is a strongly convex function with parameter $\mu>0$, Nesterov showed in \cite[Theorem 2.1.13]{Nesterov2004} that a lower bound for first-order optimisation methods of the order  $O((\frac{\sqrt{q}-1}{\sqrt{q}+1})^{2k})$ can be shown, with $q:=L/\mu\geq 1$ being the \emph{condition number} of  $F$. In this case, improved linear convergence rates of the order $O((\frac{\sqrt{q}-1}{\sqrt{q}})^{k})$ are proved. Similar results for implicit gradient descent have been studied by G\"{u}ler  \cite{Guler1992}. We also refer the reader to \cite{SalzoVilla2012}, where a general framework for inexact accelerated methods is presented.

If the objective function in \eqref{opt:problem} can be further decomposed into the sum of a convex function $f$ with Lipschitz gradient $\nabla f$ and a convex, l.s.c. and non-smooth function $g$, i.e. if the problem \eqref{opt:problem} can be rewritten as
\begin{equation}  \label{composite:opt}
\min_{x\in\mathcal{X}}~ \left\{ F(x)=f(x) + g(x)\right\} ,
\end{equation}
different descent methods taking into account the non-differentiability of $F$ need to be considered. Such approaches go under the name of \emph{composite optimisation} methods, after the work of Nesterov \cite{Nesterov2013}. A typical optimisation strategy for solving composite optimisation problems consists in alternating along the iterations a `forward' (i.e. explicit) gradient descent step taken in correspondence with the differentiable component $f$ and a `backward' (implicit)  gradient descent step in correspondence with the non-smooth part $g$. Due to this alternation, such optimisation technique is known as \emph{forward-backward} (FB) splitting. The literature on FB splitting methods is extremely vast. Historically, such strategy has firstly been used in \cite{Goldstein1964} for projected gradient descent, and subsequently popularised within the imaging community after the work of Combettes and Wajs \cite{CombettesWajs2005}. Acceleration methods for FB splitting has firstly been considered by Nesterov in \cite{Nesterov2004} for projected gradient descent, and later extended by Beck and Teboulle  \cite{BeckTeboulle2009} to more general `simple' non-smooth functions $g$ under the name of Fast Iterative Shrinkage/Thresholding Algorithm (FISTA). Several variants of FISTA have been considered in a number of work such as \cite{Nesterov2005,Tseng2008,Nesterov2013,ChambollePock2015b,Burger2016,Bonettini2016}  , just to mention a few, and properties such as convergence of the iterates under specific assumptions (\cite{ChambolleDossal2015}) and monotone variants (M-FISTA) \cite{BeckTeboulle2009b,TaoBoleyZhang2016} have also been studied. In the case when only an approximate evaluation of the FB operators up to some error can be provided, accelerated convergence rates can also be shown. We refer the reader to  \cite{SchmidtRouxBack2011,VillaSalzo2013,AujolDossal2015} for these studies

In its original formulation,  FISTA requires an estimate on the Lipschitz constant $L_f>0$ of $\nabla f$. Whenever such estimate is not easily computable, an Armijo-type backtracking rule \cite{Armijo1966} can  alternatively be used \cite[Section 4]{BeckTeboulle2009}. By construction, this backtracking strategy  requires such estimate to be non-decreasing along the iterations. From a practical point of view, this  conditions implies that if a large value of this constant is computed in the early iterations, a corresponding small (or even smaller!) gradient step size will be used in the later iterations. As a consequence, convergence speed may suffer if an inaccurate estimate of $L_f$ is computed. To avoid this drawback,  Scheinberg, Goldfarb and Bai have proposed in \cite{Scheinberg2014} a backtracking strategy for FISTA where and adaptive increasing and decreasing of the estimated Lipschitz constant along the iterations is allowed. In particular, a Lipschitz constant estimate is computed locally at each iterate $k\geq 1$ in terms of a suitable average of the $k-1$ local estimates of the $L_f$ computed in the previous iterations. The proposed strategy is shown to guarantee acceleration and to outperform the standard Armijo-type backtracking in several numerical examples. Compared to the similar full backtracking strategy proposed by Nesterov in \cite{Nesterov2013}, the criterion used in \cite{Scheinberg2014} renders cheaper since it does not require the extra calculation of the term $\nabla f$ in correspondence with the proximal step at each iteration. 

In the case of strongly convex objective functionals, improved linear convergence rates are expected. Recalling the composite problem \eqref{composite:opt}, the case of a strongly convex component $f$ has firstly been considered for projected gradient descent in \cite{Nesterov2004} and, more recently, extended by Chambolle and Pock \cite{ChambollePock2016} to the case of  strongly convex$f$ and $g$.  In this work, we will denote this general FISTA algorithm by GFISTA. For GFISTA, linear convergence rates have rigorously been shown, encompassing the quadratic ones of plain FISTA in the non-strongly convex case.  For its practical application, GFISTA requires an estimate of the Lipschitz constant $L_f$, which paves the way for the design of robust and fast backtracking strategies similar to the ones described above. We address this problem in this work.

\subsection*{Contribution} In this work we analyse a full backtracking strategy for the  strongly convex version of FISTA (GFISTA)  proposed in \cite{ChambollePock2016}. Differently from the standard backtracking rule proposed in the original paper by Beck and Teboulle \cite{BeckTeboulle2009} and based on an Armijo line-searching \cite{Armijo1966}, the strategy considered here allows for both increasing and decreasing of the Lipschitz constant estimate, i.e. for both decreasing and increasing of the gradient descent step size. Compared to the full backtracking strategy already presented by Nesterov in \cite{Nesterov2013}, the one we consider here does not require the evaluation of the gradient of the smooth component in correspondence with the proximal step at each iteration, thus it renders cheaper. A similar backtracking strategy has been considered by Scheinberg, Godfarb and Bai in \cite{Scheinberg2014} for plain FISTA, but its generalisation to the strongly convex case is not  straightforward. We address this in this work, presenting a unified framework where the standard FISTA algorithm (with and without backtracking) can be derived as a particular case. In the case of strongly convex objectives, we prove linear convergence results studying in detail the decay speed of the corresponding convergence factors. We validate our theoretical results on some exemplar problems with strongly convex objective functions which can be encountered in imaging or in data analysis. To relax the dependence on the strong convexity parameters appearing in the algorithm, we finally combine the backtracking strategy to classical restarting methods \cite{Candes2015}, which show empirical convergence properties.

\subsection*{Organisation of the paper} In Section \ref{sec:notation} we recall some definitions and standard assumptions used in the modelling of composite optimisation problems. In Section \ref{FISTA-CP} we present the GFISTA strongly convex variant of FISTA studied in \cite{ChambollePock2016}. Next, in Section \ref{sec:backtr_GFISTA} we analyse an adaptive backtracking strategy for GFISTA and prove the accelerate convergence results by means of technical tools inspired by \cite{Nesterov2004}. Numerical examples confirming our theoretical results are reported in Section \ref{sec:numerics}. In the final Section \ref{sec:conclusions} we summarise the main results of this work and give an outlook to some challenging questions to be addressed in future work.

\subsection*{Remark}

 In their recent preprint \cite{Florea2017}, Florea and Vorobyov an algorithm similar to the one described in this work as an extension of their previous work \cite{Florea2016}. The convergence result \cite[Theorem 2, Section 3.1]{Florea2017} obtained by the authors is similar to the one presented in our work (see Theorem \ref{theo:convergence}), but les accurate since it is based on a worst-case analysis, while ours depends on average quantities estimated along the iterations. Furthermore, the arguments used in  \cite{Florea2017} are completely different from the ones used here. To show the main convergence result, the authors considered  \emph{generalised estimate sequences}, a notion which, starting from the original paper by Nesterov \cite{Nesterov1983}, has indeed become very popular in the field of optimisation (see, e.g., \cite{Guler1992,LinMairal2015,SalzoVilla2012}, just to mention a few) due to its easy geometrical interpretation. However, the use of this technique leaves the technical difficulties related to the precise study of the decay speed of the convergence factors somehow hidden. Inspired by \cite{Nesterov2004} and  \cite{BeckTeboulle2009}, we follow here a different path, defining appropriate decay factors and extrapolation rules along the iterations which, eventually, will result in an accelerated (linear) convergence rates.

\section{Preliminaries and notation}  \label{sec:notation}

We are interested in the solution of the composite minimisation problem
\begin{equation}   \label{min:prob}
\min_{x\in\mathcal{X}}~ \left\{ F(x)= f(x) + g(x)\right\} ,
\end{equation}
where $\mathcal{X}$ is a (possibly infinite-dimensional) Hilbert space endowed with norm $\|\cdot\| = \langle\cdot,\cdot\rangle^{1/2}$ and $F:\mathcal{X}\to \mathbb{R}\cup\{+\infty\}$ is a convex, l.s.c. and proper functional to minimise. We denote by $x^*\in \mathcal{X}$ a minimiser of $F$. We assume that $f:\mathcal{X}\to\mathbb{R}$ is a differentiable convex function  with Lipschitz gradient and $g:\mathcal{X}\to\R\cup\left\{+\infty\right\}$ is non-smooth, convex and l.s.c .  We further denote by $L_f$ the Lipschitz constant of $\nabla f$, so that
$$
\| \nabla f(y) - \nabla f(x) \| \leq L_f \|y-x\|, \qquad\text{ for any }x,y\in\mathcal{X}.
$$
The strong convexity parameter of $f$ will be denoted by $\mu_f\geq 0$ so that for any $t\in\left[0,1\right] $, by definition, there holds
$$
f(tx+(1-t)y) \leq tf(x) + (1-t)f(y)-\frac{\mu_f}{2}t(1-t)\|x-y\|^2,\quad\text{for any }x,y\in \mathcal{X}.
$$
Similarly, by $\mu_g\geq 0$ we will denote the strong convexity parameter of $g$.
The strong convexity parameter of the composite functional $F$ in \eqref{min:prob} will be then the sum $\mu=\mu_f + \mu_g$. 

In this work we are particularly interested in the case when at least one of the two parameters $\mu_f$ and $\mu_g$ is strictly positive, so that $\mu>0$.

\subsection*{Remark}
Note that the case $\mu=0$ reduces \eqref{min:prob} to the classical FISTA-type optimisation problem. In the case of projected gradient descent, i.e. when solving
$$ 
\min_{x\in\mathcal{B}\subset\mathcal{X}}~f(x),
$$
the case $\mu_f>0$ has already been studied by Nesterov in \cite{Nesterov2004}. The problem can formulated in the form \eqref{min:prob} with $g$ being the indicator function of the subset $\mathcal{B}$ (with $\mu_g=0$) as:
$$ 
\min_{x\in\mathcal{X}}~f(x)+\delta_{\mathcal{B}}(x),\quad \text{with}\quad \delta_{\mathcal{B}}=\begin{cases}
0,\quad&\text{if }x\in\mathcal{B}\\
+\infty,\quad&\text{if }x\notin\mathcal{B}.
\end{cases}
$$
Note, however, that the proof in \cite{Nesterov2004} works actually for any function $g$, see \cite{ChambollePock2016} for more details. 

In order to write the FB optimisation step, a standard descent step in the differentiable component $f$ is combined with an implicit gradient descent step for $g$. For any $\tau >0$ and for $\bar{x}\in\mathcal{X}$ we then introduce the corresponding FB operator $T_\tau: \mathcal{X}\to\mathcal{X}$:
\begin{equation*}
\bar{x}\mapsto \hat{x} = T_{\tau}\bar{x} := \text{prox}_{\tau g}\left(\bar{x}-\tau\nabla f(\bar{x})\right),
\end{equation*}
where $\text{prox}_{\tau g}$ denotes the proximal mapping operator defined by:
$$
\text{prox}_{\tau g}(z):= \argmin_{y\in \mathcal{X} } \left( g(y) + \frac{1}{2\tau} \| z-y\|^2\right),\quad z\in\mathcal{X}.
$$
Note that in order to exploit some properties of the proximal mapping operator above, for $\eta>0$ we will also make use of the notation:
\begin{equation}  \label{prox:metric}
\text{prox}_{g}^\eta(z)= \argmin_{y\in \mathcal{X} } \left( g(y) + \frac{1}{2} \| z-y\|_{\eta^{-1}}^2\right),\quad z\in\mathcal{X},
\end{equation}
where the weighted norm is defined by $\|w\|_{\eta^{-1}}^2 = \langle \eta^{-1} w, w\rangle$.


\section{A General Fast Iterative Shrinkage/Thresholding Algorithm} \label{FISTA-CP}

The FISTA algorithm proposed in \cite{BeckTeboulle2009} is a very popular optimisation strategy to minimise composite functionals $F$ like \eqref{min:prob} with convergence guarantees of order $O(1/k^2)$. Originally proposed by Nesterov in \cite{Nesterov2004} in the case of smooth constrained minimisation, FISTA extends Nesterov's approach for more general non-smooth functions $g$. In the strongly-convex case $\mu>0$ linear convergence rates have been shown in \cite{ChambollePock2016} by means of a careful study of the decay of the composite functional towards is optimal value . In the following, we  will refer to this extension as GFISTA. 

For the sake of conciseness, we unify in Algorithm \ref{alg:GFISTA} the FISTA and GFISTA algorithms followed by the convergence result  \cite[Theorem B.10]{ChambollePock2016}. Its proof is rather technical and can be found in \cite[Appendix B]{ChambollePock2016}: the key idea consists in finding a useful recursion starting from the following descent rule for $F$ holding for every $x\in\mathcal{X}$ and for $\hat{x}=T_\tau \bar{x}$, with $\bar{x}\in\mathcal{X}$:
\begin{equation}   \label{descent:rule}
F(\hat{x}) + (1+\tau\mu_g)\frac{\|x-\hat{x}\|^2}{2\tau} \leq F(x) + (1-\tau\mu_f)\frac{\|x-\bar{x}\|^2}{2\tau}, \quad\tau>0.
\end{equation}
Inequality \eqref{descent:rule} is in fact classically used as a starting point to study convergence rates. Its proof is a trivial consequence of a general property holding for strongly convex functions. We report its proof in Lemma \ref{lemma:descent:proof} in the Appendix.

Starting from \eqref{descent:rule}, the general technique to perform a convergence analysis consists in taking as element $x\in\mathcal{X}$ the convex combination of the $k$-th iterate $x_k$ of the algorithm considered and a generic point (such as $x^*$) and, by means of (strong) convexity assumptions, in defining an appropriate decay factor by which a recurrence relation for the algorithm starting from the initial guess $x_0$ can be derived. To show acceleration, a detailed study of such factor needs then to be done by means of technical properties of the iterates of the algorithm and of its extrapolation parameters. We refer to the work of Nesterov \cite{Nesterov2004} for a review of these techniques applied to standard cases and to \cite{ChambollePock2016} to a survey on their applications in the context of Imaging.

The result reported in Theorem \ref{theo:GFISTA} generalises the ones proved for FISTA  in \cite{Nesterov2004,BeckTeboulle2009}. In particular, the standard FISTA convergence rate of $O(1/k^2)$ proved in \cite[Theorem 4.4]{BeckTeboulle2009} in the non-strongly convex case ($\mu=q=0$ and $t_0=0$) turns out to be a particular case, while improved linear convergence is shown whenever the composite functional $F$ is $\mu$-strongly convex ($\mu>0$) and an estimate on the Lipschitz constant $L_f$ is available and used as an input to find admissible gradient parameters $\tau>0$. We refer the reader to \cite{Nesterov2005,Nesterov2013,Tseng2008} for similar results proved for variants of FISTA.

\begin{algorithm}[!h]
  \caption{FISTA and GFISTA (no backtracking)}
\label{alg:GFISTA}
  \begin{algorithmic}
    \Statex {\textbf{Input}: $0<\tau\leq 1/L_f,~\mu\geq 0,~ x^0=x^{-1} \in \mathcal{X}, q :=\tau\mu/(1+\tau\mu_g)\in [0,1)$ and $t_0\in \mathbb{R}$ s.t. $0\leq t_0\leq 1/\sqrt{q}$.}
    
\vspace{0.2cm}    
    
    \For  {$k\geq 0$}
\begin{align} 
& y^k = x^k + \beta_k(x^k - x^{k-1})  \notag \\
& x^{k+1} = T_\tau y^k = \text{prox}_{\tau g} (y^k - \tau\nabla f(y^k))  \label{rules:update:xk}
\end{align} 
where:
\begin{align} \label{rules:GFISTA}
t_{k+1} = \frac{1-qt^2_k+ \sqrt{(1-qt_k^2)^2+4 t^2_k}}{2}  \\
\beta_k = \frac{t_k-1}{t_{k+1}} \frac{1+\tau\mu_g - t_{k+1}\tau\mu}{1-\tau\mu_f} \notag
\end{align}
    \EndFor 
  \end{algorithmic}
\end{algorithm}

\subsection{Remark (FISTA updates)} \label{remark:alg:FISTA}
Note that in the case $\mu=0$ the update rules for $t_{k+1}$ and $\beta_k$ in \eqref{rules:GFISTA} simplify to:
\begin{equation}   \label{rules:FISTA}
t_{k+1} = \frac{1+ \sqrt{1+4 t^2_k}}{2}, \qquad \beta_k = \frac{t_k-1}{t_{k+1}},
\end{equation}
which are the standard FISTA updates considered by Beck and Teboulle in \cite{BeckTeboulle2009}.

\newtheorem{convergence_results}{Theorem}[section]

\begin{convergence_results}[\cite{Nesterov2004} and Theorem B.1 \cite{ChambollePock2016}] \label{theo:GFISTA}
Let $\tau>0$ with $\tau \leq 1/L_f$  and let $q:=\frac{\mu\tau}{1+\tau\mu_g}$ and $x^*$ be a minimiser of $F$. If $\sqrt{q}t_0\leq 1$ with $t_0\geq 0$, then the sequence $( x^k)$ produced by the Algorithm \ref{alg:GFISTA} in \eqref{rules:update:xk} satisfies
\begin{equation*}
F(x^k)-F(x^*) \leq r_k(q) \left( t^2_0(F(x^0) - F(x^*)) + \frac{1+\tau\mu_g}{2} \| x-x^*\|^2\right),
\end{equation*}
and $r_k(q)$ is defined by:
\begin{equation*}
r_k(q) = \min\left\{\frac{4}{(k+1)^2},(1+\sqrt{q})(1-\sqrt{q})^k, \frac{(1-\sqrt{q})^k}{t^2_0}\right\}.
\end{equation*}
\end{convergence_results}

\subsection*{Backtracking} Whenever an estimate of $L_f$ is not available, backtracking techniques can be used.  For FISTA, an Armijo-type backtracking rule has been proposed in the original paper of Beck and Teboulle \cite{BeckTeboulle2009}. For that, similar convergence rates as above can be proved. Furthermore, in order to improve the speed of the algorithm allowing also the increasing of the step size $\tau$ in the neighbourhoods of `flat' points of the function $f$ (i.e. where $L_f$ is small), a full backtracking strategy for FISTA has been considered by Scheinberg, Goldfarb and Bai in \cite{Scheinberg2014}. 

The typical inequality to check in the design of any backtracking strategy can be derived from \eqref{descent:rule} (see Lemma \ref{lemma:descent:proof} in the Appendix) and reads:

\begin{equation}   \label{condition:backtr}
F(\hat{x}) + (1+\tau\mu_g)\frac{\|x-\hat{x}\|^2}{2\tau} + \left(  \frac{\|\hat{x}-\bar{x}\|^2}{2\tau} - D_f(\hat{x},\bar{x}) \right) \leq F(x) + (1-\tau\mu_f)\frac{\|x-\bar{x}\|^2}{2\tau},
\end{equation}

\noindent where $D_f(\hat{x},\bar{x}):=f(\hat{x})-f(\bar{x})-\langle \nabla f(\bar{x}),  \hat{x}-\bar{x} \rangle \leq \frac{L_f}{2} \|\hat{x}-\bar{x}\|^2$ is the Bregman distance of $f$ between $\hat{x}$ and $\bar{x}$. Note that in the case when no backtracking is performed, condition \eqref{condition:backtr} is satisfied as long as:
\begin{equation}   \label{condition:bregman}
D_f(\hat{x},\bar{x})\leq \frac{\|\hat{x}-\bar{x}\|^2}{2\tau},    \tag{CB}
\end{equation}
which is clearly true for constant $\tau$ whenever $0<\tau\leq 1/L_f$ with $L_f$ known. However, by letting $\tau$ vary, one can alternatively check condition \eqref{condition:backtr} along the iterations of the algorithm and redefine $\tau_k$ at each iteration $k\geq 1$ so as to compute a local Lipschitz constant estimate. 

In the following, we will indeed use this rule for the design of a backtracking strategy for Algorithm \ref{alg:GFISTA} with $\mu>0$. In order to allow robust backtracking, we will allow the step size $\tau_k$ to either decrease (as it is classically done) or increase depending on the validity of the following inequality:
\begin{equation}   \label{condition:increase}
\frac{2D_f(\hat{x},\bar{x})}{\|\hat{x}-\bar{x}\|^2}>\rho\left(\frac{1}{\tau_k}\right),  \tag{CB2}
\end{equation}
where the constant $\rho\in (0,1)$ is chosen in advance. Note that this inequality entails that at any iteration the following inequality holds:
\begin{equation}    \label{eq:Lk_rho}
\tau_k \geq \frac{\rho}{L_f}.
\end{equation}
Heuristically, condition \eqref{condition:increase} favours the step size $\tau_k$ to be decreased at iteration $k\geq 1$ whenever the estimate of the Lipschitz constant given by the left hand side in the inequality above is `too close' to $1/\tau_k$, i.e. whenever \eqref{condition:increase} is verified, and increased otherwise. 

\section{A backtracking strategy for GFISTA algorithm \ref{alg:GFISTA}} \label{sec:backtr_GFISTA}

Following the analysis performed in \cite[Section 4, Appendix B]{ChambollePock2016}, we prove that the backtracking strategy described above and applied to the GFISTA algorithm \ref{alg:GFISTA} enjoys accelerated convergence rates, which turn out to be linear in the case $\mu>0$.


For an arbitrary $t\geq 1$, $k\geq 0$ and $\tau>0$ we start from inequality \eqref{condition:backtr} and choose the point $x$ to be the convex combination $x=((t-1)x^{k} + x^*)/t$ where  $x^{k}$ is an iterate of the algorithm we are going to define and $x^*$ is a minimiser of $F$. For the other points, we set $\bar{x}=y^{k+1}$ and $\hat{x}=x^{k+1}=T_\tau y^{k+1}$. The formula for $y^{k+1}$ will be specified in the following. 

After multiplication by $t^2$ and using the strong convexity of $F$ we get:
\begin{multline} \label{decay1}
t^2\left( F(x^{k+1})-F(x^*)\right) + \frac{1+\tau\mu_g}{2\tau} \|x^*-x^{k+1}-(t-1)(x^{k+1}-x^k)\|^2\\ + ~  t^2 (t-1)\frac{\mu(1-\tau\mu_f)}{1+\tau\mu_g-t\tau\mu}\frac{\|x^k-y^{k+1}\|^2}{2} \leq  t(t-1)\left( F(x^k) - F(x^*) \right) \\ +~ \frac{1+\tau\mu_g-t\tau\mu}{2\tau} \| x^*-x^k - t\frac{1-\tau\mu_f}{1+\tau\mu_g-t\tau\mu}(y^{k+1}-x^k) \|^2.
\end{multline}
We now set $t=t_{k+1}$, let $\tau = \tau_{k+1}$ and define the following quantities:
\begin{align}
\tau'_{k+1} &:= \frac{\tau_{k+1}}{1+\tau_{k+1}\mu_g} >0 \label{def_tauk}\\
q_{k+1} &:= \mu\tau'_{k+1} = 1 - \frac{1-\tau_{k+1}\mu_f}{1+\tau_{k+1}\mu_g}\in [0,1),   \label{def:q}\\
\omega_{k+1} & := \frac{1+\tau_{k+1}\mu_g-t_{k+1}\tau_{k+1}\mu}{1+\tau_{k+1}\mu_g} = 1 - t_{k+1}q_{k+1}\in (0,1],  \label{def:omega} \\
\beta_{k+1} & := \frac{t_k-1}{t_{k+1}} \frac{1+\tau_{k+1}\mu_g - t_{k+1}\tau_{k+1}\mu}{1-\tau_{k+1}\mu_f} = \omega_{k+1}\frac{t_k -1}{t_{k+1}} \frac{1+\tau_{k+1}\mu_g}{1-\tau_{k+1}\mu_f}
, \label{def:beta}
\end{align}
where we can assume $\mu_f<L_f$, so that $\tau<1/L_f$. 

We now define the following update for $y^{k+1}$:
\begin{equation}  \label{condition_y}
y^{k+1} = x^k + \beta_{k+1}(x^k-x^{k-1}), 
\end{equation}
for any $k\geq 0$.  After further multiplying  \eqref{decay1} by $\tau'_{k+1}$, we thus  deduce:
%
\begin{multline}   \label{eq:decay_new}
\tau_{k+1}'t_{k+1}^2\left( F(x^{k+1})-F(x^*)\right) + \frac{1}{2}\|x^*-x^{k+1}-(t_{k+1}-1)(x^{k+1}-x^k)\|^2\\ \leq \tau_{k+1}'t_{k+1}(t_{k+1}-1)\left( F(x^k) - F(x^*) \right)  \\ +~ \frac{\omega_{k+1}}{2} \| x^*-x^{k}-(t_{k}-1)(x^{k}-x^{k-1})  \|^2.
\end{multline}

%

Let us now assume that for every $k\geq 1$ the following inequality holds:
\begin{equation}  \label{inequalityFISTA}
\tau'_{k+1}t_{k+1}(t_{k+1}-1)\leq \omega_{k+1}\tau'_k t^2_k,
\end{equation}
and that the same holds for the iteration $k=0$ by defining $T_0^2:=\tau_0' t_0^2$ implicitly by
\begin{equation}  \label{equality:first:iter}
T_0^2 = \frac{\tau'_1 t_1(t_1 -1)}{\omega_1} =  \frac{\tau_1 t_1(t_1 -1)}{1+\tau_1\mu_g - t_1\tau_1\mu},
\end{equation}
which is positive whenever
\begin{equation}   \label{cond:t1}
1\leq t_1 < \frac{1+\tau_1\mu_g}{\tau_1\mu} = \frac{1}{q_1}.
\end{equation}
Then, we get from \eqref{eq:decay_new} that for any $k\geq 0$:
\begin{multline} \label{decay4}
\tau_{k+1}'t_{k+1}^2\left( F(x^{k+1})-F(x^*)\right) + \frac{1}{2}\|x^*-x^{k+1}-(t_{k+1}-1)(x^{k+1}-x^k)\|^2\\ \leq \omega_{k+1}\left(\tau_k't^2_{k}\left( F(x^k) - F(x^*) \right) +~ \frac{1}{2} \| x^*-x^{k}-(t_{k}-1)(x^{k}-x^{k-1})  \|^2\right).
\end{multline}

By now applying \eqref{decay4} recursively and if we let $x^0 = x^{-1}\in\mathcal{X}$, we  find the following convergence inequality 
\begin{equation} \label{decay5}
F(x^{k})-F(x^*) \leq \theta_k~\left(T_0^2\left( F(x^0) - F(x^*) \right) +~ \frac{1}{2} \| x^*-x^0 \|^2\right),
\end{equation}
where the decay rate of the factor 
\begin{equation}   \label{def:theta_k}
\theta_k : = \frac{\displaystyle\prod_{i=1}^{k} \omega_{i}}{\tau'_kt^2_k}
\end{equation}
needs to be studied to determine the speed of convergence of $F(x^k)$ to the optimal value $F(x^*)$. We will do this in the following sections using some technical properties of the sequences defined above.

\subsection{Update rule} 

Assuming that \eqref{inequalityFISTA} holds with an equality sign, i.e. if
\begin{equation}  \label{equality_tk}
\tau'_{k+1}t_{k+1}(t_{k+1}-1) = \omega_{k+1}\tau'_k t^2_k,
\end{equation}
and after recalling the definition of $\omega_{k+1}$ in \eqref{def:omega},  we find the following update rule for the elements of sequence $(t_k)$, $k\geq 1$:
\begin{align}  
t_{k+1} & = \frac{1-q_{k+1}\frac{\tau'_{k}}{\tau'_{k+1}}t^2_k + \sqrt{\left(1-q_{k+1}\frac{\tau'_{k}}{\tau'_{k+1}}t^2_k \right)^2 + 4\frac{\tau'_{k}}{\tau'_{k+1}}t^2_k}} {2}  \notag \\
& = \frac{1-q_k t^2_k + \sqrt{\left(1-q_k t^2_k \right)^2 + 4\frac{q_{k}}{q_{k+1}}t^2_k}} {2} \geq 0, \label{updatetk}
\end{align}
by \eqref{def:q} and \eqref{def_tauk}.

We can now present the  GFISTA algorithm with backtracking.

\begin{algorithm}[H]
  \caption{GFISTA with backtracking}
\label{alg:GFISTA:backtr}
  \begin{algorithmic}
    \Statex {\textbf{Input}:  $\mu_f$, $\mu_g$, $\tau^0_1>0$, $q_1:=\mu \tau^0_1/(1+\tau^0_1\mu_g)$, $\rho\in(0,1)$, $y^1=x^0=x^{-1} \in \mathcal{X}$ and $t_1\in\R$ s.t. $1\leq t_1\leq 1/\sqrt{q_1}$.}
    
    \vspace{0.2cm}

    \For  {$k\geq 1$}

Set $i=0$;
\If { condition \eqref{condition:increase} is satisfied}
\While {condition \eqref{condition:bregman} is \textbf{not} satisfied \textbf{and} $i\leq i_{max}$}

{
\qquad\quad \textbf{reduce step-size}: $\tau_{k} = \rho^{i}~\tau^0_k$;

\qquad\quad Compute 
\begin{equation} \label{rules:update:xk1}
x^{k} = T_{\tau_{k}}~ y^{k} = \text{prox}_{\tau_{k} g} (y^{k} - \tau_{k}\nabla f(y^{k})) 
\end{equation}

\qquad\quad $i=i+1$;

}

\EndWhile
\ElsIf{\eqref{condition:increase}  is \textbf{not} satisfied}

{\textbf{increase step-size}: $\tau_{k}=\frac{\tau^0_{k}}{\rho}$ ;

Compute $x_{k}$ using \eqref{rules:update:xk1};

}

\EndIf

\noindent Set
$$
\tau_{k+1}^0= \tau_k,\qquad \qquad q_{k+1}=\frac{\mu\tau^0_{k+1}}{1+\tau^0_{k+1}\mu_g}.
$$
Compute $t_{k+1}$ using  \eqref{updatetk}.

\noindent Set
\begin{align}   \label{updates:y:b}
\beta_{k+1} & = \frac{t_k-1}{t_{k+1}} \frac{1+\tau^0_{k+1}\mu_g - t_{k+1}\tau^0_{k+1}\mu}{1-\tau^0_{k+1}\mu_f}.\\
 y^{k+1} & = x^{k} + \beta_{k+1}(x^{k}-x^{k-1}). 
\end{align} 
    \EndFor 
  \end{algorithmic}
\end{algorithm}

We remark that compared to the algorithm studied in \cite[Section 4]{Nesterov2013}, Algorithm \ref{alg:GFISTA:backtr} has a lower per-iteration cost. The reason for that is that the backtracking criterion considered in \cite{Nesterov2013} requires at any iteration $k$ the computation of the quantity $\nabla f (T_{\tau_{k+1}} y^k)$, whereas our backtracking condition \eqref{condition:bregman} is based on the calculation of $D_f$, and the sole computation of $\nabla f(y^k)$ is required, thus avoiding the calculation of $\nabla f$ in the proximal step. In many applications (e.g. compressed sensing), this difference can be quite crucial: the extra-evaluation of $\nabla f$ in one point requires in fact two matrix-vector multiplications compared to a single one required for functional evaluation.
Similar considerations have already been made for the FISTA algorithm with full backtracking in \cite{Scheinberg2014} since the stopping criterion for the backtracking procedure considered therein is in fact similar to the one used in our Algorithm \ref{alg:GFISTA:backtr}.

\subsection{Remark (No backtracking)} \label{remark:update_ruleFISTA}
When no backtracking is performed along the iterations $\tau_k = \tau_{k+1}$ for any $k$ and the ratio $q_k/q_{k+1}$ in \eqref{updatetk} is constantly equal to one. In this case, the update rule \eqref{updatetk} is the same as the one used in \eqref{rules:GFISTA} for GFISTA without backtracking, compare \cite[Appendix B]{ChambollePock2016} . 

In the non-strongly convex case ($q_{k}=0$ for every $k$), the update rule \eqref{updatetk} is exactly the same \eqref{rules:FISTA} for the original FISTA algorithm  \cite{Nesterov2004,BeckTeboulle2009}.

\subsection{Remark (FISTA with backtracking)}
In the non-strongly convex case ($\mu_f=\mu_g=q_{k}=0$ for every $k$), \eqref{updatetk} reduces to
$$
t_{k+1} = \frac{1 + \sqrt{1 + 4  \frac{\tau_{k}}{\tau_{k+1}} t^2_k} }{2},
$$
which is exactly the same update rule considered by Goldfarb et al. in \cite{Scheinberg2014} for adaptive backtracking of plain FISTA.


We now prove a fundamental property of the sequence $(t_k)$ defined by \eqref{updatetk}.

\newtheorem{tkgreater}{Lemma}[section]

\begin{tkgreater}   \label{tkgreaterzero}
Let the sequence $(t_k)$ be defined by the update rule \eqref{updatetk}. Then: 
$$
t_k\geq 1\qquad \text{for any }k\geq 1.
$$
\end{tkgreater}

\begin{proof}
We simply observe that since $q_{k}\leq 1$ for every $k$ there holds:
\begin{align}  
t_{k} &= \frac{1-q_{k-1}t^2_{k-1} + \sqrt{\left(1-q_{k-1}t^2_{k-1} \right)^2 + 4\frac{q_{k-1}}{q_{k}}t^2_{k-1}}} {2} . \notag\\
&\geq \frac{1-q_{k-1}t^2_{k-1} + \sqrt{\left(1-q_{k-1}t^2_{k-1} \right)^2 + 4 q_{k-1}t^2_{k-1}}} {2} . \notag\\
& =  \frac{1-q_{k-1}t^2_{k-1} + \sqrt{\left(1+q_{k-1}t^2_{k-1}\right)^2}} {2}=1. \notag
\end{align}
\end{proof}

%
%
%

For the following convergence proofs, the following technical lemma will be crucial.

\newtheorem{propertyqktk}{Lemma}[section]

\begin{propertyqktk}   \label{lemma:propertyqktk}
Let $\sqrt{q_1}t_1\leq 1$. Then, there holds:
\begin{equation}  \label{lemma_qktk}
\sqrt{q_k}t_k \leq 1 .
\end{equation}
\end{propertyqktk}

\begin{proof}


We proceed by induction. By assumption, the initial step $k=1$ holds. Let us assume that \eqref{lemma_qktk} holds for some $k\geq 1$. By \eqref{equality_tk}, we get:
$$
q_{k+1} t^2_{k+1} = q_{k+1} t_{k+1} + \omega_{k+1}q_kt^2_k = 1+ \omega_{k+1}(q_kt^2_k -1 ) \leq 1
$$
by simply applying the induction assumption.
\end{proof}

Note that the condition $t_1 \leq 1/\sqrt{q_1}$ combined with $t_1\geq 1$  results in the following bound:
\begin{equation}  \label{bounds:t1}
1\leq t_1 \leq \sqrt{1 + \frac{1-\tau_1\mu_f}{\tau_1\mu}}.
\end{equation}
Furthermore, since $1/\sqrt{q_1} < 1/q_1$, such condition also guarantees  \eqref{cond:t1}. In particular, $t_1=1$ is an admissible choice.

\subsection{Convergence rates} 
In this section, we follow \cite{Nesterov2004,ChambollePock2016} to derive a precise estimate of the factor $\theta_k$ in \eqref{def:theta_k}. 


The following convergence result shows that the backtracking strategy applied to the GFISTA algorithm guarantees accelerated linear convergence rates given in terms of \emph{averaging} quantities defined in terms of the Lipschitz constant estimates along the iterations. Comments on our result in comparison to the ones studied in analogous works \cite{Scheinberg2014,Florea2017} are given in the following remarks.

\newtheorem{convergence_rates}{Theorem}[section]

\begin{convergence_rates}[Convergence rates]   \label{theo:convergence}
Let $T_0$ be defined as in \eqref{equality:first:iter}. If $1\leq t_1 \leq 1/\sqrt{q_1}$, then the sequence $(x^k)$ produced by the  Algorithm \ref{alg:GFISTA:backtr} with \eqref{updatetk}, \eqref{def:omega}, \eqref{def:beta} and \eqref{condition_y} satisfies:
\begin{equation}  \label{convergence:functional2}
F(x^k)-F(x^*)\leq r_k \left( T_0^2 \left(F(x^0)-F(x^*)\right) + \frac{1}{2}\|x^0-x^*\|^2\right),
\end{equation}
where $r_k$ is defined by:
\begin{equation}   \label{convergence:rate}
r_k := \min\left\{\frac{4 \bar{L}_k}{k^2}, (L_1 - \mu_f)(1-\sqrt{\bar{q}_k})^{k-1} \right\},
\end{equation}
and the average quantities $\bar{L}_k$ and $\sqrt{\bar{q}_k}$ are  defined by:
\begin{equation}   \label{eq:Lbar}
\sqrt{\bar{L}_k} := \frac{1}{\frac{1}{k} \sum_{i=1}^k  \frac{1}{\sqrt{L_i - \mu_f}}},\qquad \sqrt{ \bar{q}_k }:= \frac{1}{k-1}\sum_{i=2}^{k} \sqrt{ \frac{\mu}{L_i + \mu_g} },
\end{equation}
with $L_i:=1/\tau_i$.

\end{convergence_rates}

\begin{proof}
We recall the definition of $\theta_k$ given in \eqref{def:theta_k} and start computing the $O(1/k^2)$ factor in \eqref{convergence:rate} following \cite{Nesterov2004,ChambollePock2016}.

We first notice that from \eqref{equality_tk} we can deduce
\begin{equation}  \label{relation:omega}
1-\frac{1}{t_{k+1}} = \omega_{k+1}\frac{\tau'_{k}t^2_{k}}{\tau'_{k+1}t^2_{k+1}} = \frac{\theta_{k+1}}{\theta_k} \leq 1,
\end{equation}
which also shows that $\theta_k$ is non-increasing.
Thus, we have:
\begin{equation}   \label{eq:estim1}
\frac{1}{\sqrt{\theta_{k+1}}} - \frac{1}{\sqrt{\theta_{k}}} = \frac{\theta_k-\theta_{k+1}}{\sqrt{\theta_k\theta{k+1}}(\sqrt{\theta_k} + \sqrt{\theta_{k+1})}} \geq \frac{\theta_k - \theta_{k+1}}{2\theta_k\sqrt{\theta_{k+1}}}.
\end{equation}
By now applying \eqref{relation:omega}, we get
$$
\frac{1}{\sqrt{\theta_{k+1}}} - \frac{1}{\sqrt{\theta_{k}}} \geq \frac{1}{2 t_{k+1}\sqrt{\theta_{k+1}}}.
$$
We now recall definitions \eqref{def:omega}, \eqref{def:theta_k}, and use Lemma \ref{tkgreaterzero} to find:
\begin{align}
t_{k+1}\sqrt{\theta_{k+1}} & = \frac{1}{\sqrt{\tau'_{k+1}}}\prod_{i=1}^{k+1} \sqrt{\omega_i} \leq \sqrt{ \frac{\omega_{k+1}}{\tau'_{k+1}} } =  \sqrt{  \frac{1}{\tau'_{k+1}}-\mu t_{k+1} } \notag \\
& \leq \sqrt{  \frac{1}{\tau'_{k+1}}-\mu } = \sqrt{\frac{1}{\tau_{k+1} } -\mu_f},  \notag
\end{align}
whence:
$$
\frac{1}{\sqrt{\theta_{k+1}}} - \frac{1}{\sqrt{\theta_{k}}}  \geq \frac{1}{2\sqrt{\frac{1}{\tau_{k+1}}- \mu_f }}.
$$
Applying this recursively we get that for any $k\geq 1$
\begin{equation}  \label{sum:theta}
\frac{1}{\sqrt{\theta_k}} \geq \frac{1}{2}\sum_{i=1}^k \frac{1}{\sqrt{\frac{1}{\tau_i} -\mu_f}}.
\end{equation}
Note that indeed for $i=1$ we have:
\begin{equation}   \label{boundtheta1}
\theta_1 = \frac{1-\mu t_1\tau_1'}{\tau_1't^2_1} = \frac{1-\mu_g(t_1 -1)\tau_1-\mu_f t_1\tau_1}{\tau_1 t_1^2}\leq \frac{1}{\tau_1}-\mu_f ,
\end{equation}
since $t_1\geq 1$ by \eqref{bounds:t1}. We then deduce:
\begin{equation*}
\frac{1}{\sqrt{\theta_1}} \geq \frac{1}{2\sqrt{\frac{1}{\tau_1}-\mu_f}}.
\end{equation*}
After setting $L_i = 1/\tau_i$ in \eqref{sum:theta}, we get:
\begin{equation}   \label{eq:estimm2}
\sqrt{\theta_k} \leq \frac{2}{k} \sqrt{\bar{L}_k}
\end{equation}
where  $\sqrt{\bar{L}_k}$ is defined in \eqref{eq:Lbar}. 

To get the linear rates, we notice that by Lemma \ref{lemma:propertyqktk}, relation \eqref{relation:omega} and definition  \eqref{def:q}, we have:
\begin{align}   \label{eq:thetak}
 \theta_k & = \theta_1 \prod_{i=2}^k  \left(1-\frac{1}{t_i} \right)   \leq \theta_1\prod_{i=2}^k \left(1-\sqrt{q_i} \right) \\ 
 & \leq \theta_1 \prod_{i=2}^k \left(1-\sqrt{ \frac{\mu}{L_i+\mu_g} }~ \right) 
\leq (\theta_1 (1-\sqrt{\bar q})^{k-1},
\end{align}
where $\sqrt{\bar{q}_k}$ is defined as in  \eqref{eq:Lbar}. and by the concavity of the function logarithm. We then get from \eqref{eq:thetak} that:
$$
\theta_k \leq \theta_1 (1-\sqrt{\bar{q}_k})^{k-1} \leq (L_1 -\mu_f) (1-\sqrt{\bar{q}_k})^{k-1}
$$
by \eqref{boundtheta1}.  Combining this with \eqref{eq:estimm2} we finally get the final rate \eqref{convergence:rate}.
\end{proof}
%


Note that the averaging term $\bar{L}_k$ appearing above is always smaller than the actual average of the terms $(L_i-\mu_f)$, since:
\begin{equation}   \label{ineq:mean}
 \sqrt{\bar{L}_k}  \leq \frac{1}{k} \sum_{i=1}^k \sqrt{L_i - \mu_f}  \leq \sqrt{ \frac{1}{k}\sum_{i=1}^k (L_i - \mu_f)  }.
\end{equation}
Furthermore, whenever $L_f$ is known and recalling \eqref{eq:Lk_rho}, we can deduce the following bounds for the terms defined in \eqref{eq:Lbar}:
$$
\sqrt{\bar{L}_k}\leq \sqrt{\frac{L_f -\rho \mu_f}{\rho}},\qquad\qquad \sqrt{\bar{q}_k}\geq \sqrt{\frac{\rho\mu}{L_f + \rho\mu_g}}.
$$
Hence, the convergence rate $r_k$ in\eqref{convergence:rate} can be estimated as:
$$
r_k \leq \frac{1}{\rho}\min\left\{\frac{4(L_f-\rho\mu_f)}{ k^2},(L_f-\rho\mu_f)\left(1-  \sqrt{\frac{\rho\mu}{L_f + \rho\mu_g}}\right)^{k-1}  \right\}.
$$
Finally, as far as the choice of $T_0$ is concerned, note that by  \eqref{equality:first:iter} when $t_1=1$, then $T_0=0$.


\subsection*{Remark (FISTA with backtracking)}
Note that in the non-strongly convex case ($\mu=q_k=0$ for all $k$), the global convergence rate \eqref{convergence:functional2}-\eqref{convergence:rate} is analogous to \cite[Theorem 3.3]{Scheinberg2014}, which reads:
\begin{equation}   \label{godfarb:estimate}
F(x^k)-F(x^*)\leq \frac{2 \tilde{L}_k \|x^0-x^*\|^2}{\rho k^2}, 
\end{equation}
and where the term $\tilde{L}_k$ is defined by
$$
\tilde{L}_k:= \frac{\left(\sum_{i=1}^k \sqrt{L_i}\right)^2  }{k^2}.
$$
Note in fact that whenever $\mu_f=0$ our definition \eqref{eq:Lbar} relates with the one above via Remark \eqref{ineq:mean}.
 

\subsection*{Remark}
The worst-case convergence result \cite[Theorem 2]{Florea2017} is obtained via the analysis of generalised estimate sequences. In \cite[Section 4]{Florea2017} some comments on the extrapolated form of their algorithm and its relation with the strongly-convex variant of the FISTA algorithm \ref{alg:GFISTA} are given. Although the expression of the sequence $\left\{\omega_k\right\}$ and the update rule for the elements $\left\{t_k\right\}$ is similar (but not equal) to our definitions \eqref{def:omega} and \eqref{updatetk}, respectively, the arguments used by the authors are different from the ones used here. More importantly, compared to a worst-case analysis, the convergence result \ref{theo:convergence} is more precise, since it provides quantitative convergence estimates in terms of the average quantities $\sqrt{\bar{L}_k}$ and $\sqrt{\bar{q}_k}$ estimated along the iterations.

\subsection{Monotone algorithms} \label{sec:monotone}

As already noticed for standard FISTA \cite[Section V.A]{BeckTeboulle2009b} and for GFISTA without backtracking \cite[Remark B.3]{ChambollePock2016}, the convergence of the composite energy $F$ to the optimal value $x^*$ is not guaranteed to be monotone non-increasing. A straightforward modification of the GFISTA Algorithm \ref{alg:GFISTA:backtr} enforcing such property and used in several papers \cite{BeckTeboulle2009b,Tseng2008} consists in taking as $x^{k}$ any point such that $F(x^{k})\leq F(T_{\tau_{k}} y^k)$.
Recalling the definition of $\omega_{k+1}$ in \eqref{def:omega}, the update rule \eqref{updates:y:b} for extrapolation can then be changed as:
\begin{align}  \label{monotone:update}
y^{k+1} & =x^k + \beta_{k+1} \left(x^k - x^{k-1}\right) + \omega_{k+1}\frac{t_k}{t_{k+1}}\frac{1+\tau_{k+1}\mu_g}{1-\tau_{k+1}\mu_f}\left( T_{\tau_{k}} y^{k} - x^k\right)  \tag{C$2_{m}$}\\
& = x^k + \beta_{k+1}\left( \left(x^k - x^{k-1}\right) + \frac{t_{k}}{t_k-1}\left( T_{\tau_{k}} y^{k} - x^k\right)\right). \notag
\end{align}
One can easily check that starting from \eqref{descent:rule} and replacing in \eqref{decay1}  $x^{k+1}$ by $T_{\tau} y^{k+1}$ with the update rule above  the same computations of the previous sections carry on and the same convergence rates are obtained. Condition \eqref{monotone:update} suggests also a natural choice for $x^{k}$. In fact, one can simply set:
\begin{equation}  \label{monotone:update2}
x^{k} = \begin{cases}
T_{\tau_{k}} (y^{k}) \qquad & \text{if } F(T_{\tau_{k}} y^{k})\leq F(x^{k-1}),\\
x^{k-1}\qquad & \text{otherwise},
\end{cases}
\end{equation}
so that in either case one of the two terms in \eqref{monotone:update} vanishes. Whenever the evaluation of the composite functional $F$ is cheap, this choice seems to be the most sensible. Another monotone implementation of FISTA has been recently considered in \cite{TaoBoleyZhang2016} where despite the further computational costs required to compute the value $x^{k}$, an empirical linear convergence rate is observed also for standard FISTA applied to strongly convex objectives. A rigorous proof of such convergence property is an interesting question of future research.

\section{Numerical examples}  \label{sec:numerics}

In this section we report some numerical experiments to confirm numerically the convergence result \ref{theo:convergence} of Algorithm \ref{alg:GFISTA:backtr}. We also discuss some heuristic restarting strategies \cite{Candes2015} in the case when the strong convexity parameters are unknown.

\subsection{TV-Huber ROF denoising} \label{sec:gauss:den}

We start considering a strongly convex variant of the well-know Rudin, Osher and Fatemi image denoising model \cite{rudinosherfatemi}  based on the use of Total Variation (TV) regularisation. In its discretised form and for a given noisy image $u^0\in\R^{m\times n}$ corrupted by Gaussian noise with zero mean and variance $\sigma^2$, the original ROF model reads:
\begin{equation}  \label{ROF:model}
\min_u~ \lambda \| Du \|_{p,1} + \frac{1}{2}\|u-u^0\|_2^2.
\end{equation}
Here, $Du = ( (Du)_1, (Du)_2)$ is the gradient operator discretised using forward finite differences (see, e.g., \cite{chambolle2004algorithm}) and the discrete TV regularisation is defined by:
\begin{equation}   \label{discrete:TV}
\| Du \|_{p,1} = \sum_{i=1}^m \sum_{j=1}^n |(Du)_{i,j}|_p =  \sum_{i=1}^m \sum_{j=1}^n \left( (Du)_{i,j,1}^p + (Du)_{i,j,2}^p \right)^{1/p},
\end{equation}
where the value of the parameter $p$ allows for both anisotropic ($p=1$) and isotropic ($p=2$) TV, which is generally preferred to reduce grid bias. The regularisation parameter $\lambda>0$ in \eqref{ROF:model} weights the action of TV-regularisation against the fitting with the Gaussian data given by the $\ell^2$ squared term.

Taking $p=2$ in \eqref{discrete:TV}, we now follow \cite[Examples 4.7 and 4.14]{ChambollePock2016} and consider a similar denoising model where a strongly convex variant of TV is employed. This can be obtained, for instance, using the $C^1$-Huber smoothing function $h_\varepsilon: \R\to \R$ defined for a parameter $\varepsilon>0$ by:
$$
h_{\varepsilon}(t):=\begin{cases}
\frac{t^2}{2\varepsilon}\qquad &\text{for }|t|\leq\varepsilon, \\
|t|-\frac{\varepsilon}{2}\qquad &\text{for }|t|>\varepsilon.
\end{cases}
$$
Applying such smoothing to the TV energy \eqref{discrete:TV} removes the singularity in a neighbourhood zero by means of a quadratic term and leaves the TV term almost unchanged otherwise. The resulting Huber-ROF image denoising model then reads:
\begin{equation*} 
\min_u ~\lambda H_\varepsilon(u)+\frac{1}{2}\|u-u^0\|_2^2,
\end{equation*}
with
\begin{equation}  \label{Huber:ROF:regulariser}
H_\varepsilon(u) :=  \sum_{i=1}^m \sum_{j=1}^n h_\varepsilon\left( \sqrt{(Du)_{i,j,1}^2 + (Du)_{i,j,2}^2 }\right).
\end{equation}

The dual problem of \eqref{Huber:ROF:regulariser} reads:
\begin{equation}  \label{Huber:ROF:dual}
\min_{\bm{p}} ~\frac{1}{2}\|D^*\bm{p} - u^0 \|^2_2 + \frac{\varepsilon}{2\lambda}\|\bm{p}\|_2^2 + \delta_{  \left\{\|\cdot\|_{2,\infty} \leq \lambda\right\} }(\bm{p}),
\end{equation}
where $\bm{p}$ is the dual variable, $D^*$ is the adjoint operator of $D$ (i.e. the discretised negative finite-difference divergence operator) and $\delta_{  \left\{\|\cdot\|_{2,\infty} \leq \lambda\right\} }$ is the indicator function defined by:
$$
\delta_{  \left\{\|\cdot\|_{2,\infty} \leq \lambda\right\} }(\bm{p}) = \begin{cases}
0\qquad & \text{if }|\bm{p}_{i,j}|_2\leq \lambda\text{ for any }i,j,\\
+\infty \qquad & \text{ otherwise}.
\end{cases}
$$
Note that \eqref{Huber:ROF:dual} is the sum of a function $f$ with Lipschitz gradient and a non-smooth function $g$ which are respectively given by:
$$
f(\bm{p})=\frac{1}{2}\|D^*\bm{p} - u^0\|^2_2, \qquad g(\bm{p})= \frac{\varepsilon}{2\lambda}\|\bm{p}\|_2^2 + \delta_{  \left\{\|\cdot\|_{2,\infty}\leq \lambda\right\} }(\bm{p}).
$$
The gradient of the differentiable component $f$ reads:
$$
\nabla f(\bm{p}) = D(D^*\bm{p}-u^0),
$$
and it is easy to show that its Lipschitz constant $L_f$ can be estimated as $L_f\leq 8$, see, e.g. \cite{chambolle2004algorithm}. Note also that $\mu_f=0$.

The function $g$ is strongly convex with parameter $\mu_g=\mu=\varepsilon/\lambda$ and its proximal map $\hat{\bm{p}} = \text{prox}_{\tau g}(\tilde{\bm{p}})$ can be easily computed pixel-wise as:
$$
\hat{\bm{p}}_{i,j}=\frac{(1+\tau\mu_g)^{-1}\bm{\tilde{p}}_{i,j}}{\max\left\{1,(\lambda(1+\tau\mu_g))^{-1}|\bm{\tilde{p}}_{i,j}|_2\right\}}, \quad \text{ for any }i, j,
$$
since, due general properties of proximal maps with added squared $\ell^2$ terms (see Lemma \ref{lemma:proximal} in the Appendix), there holds:
$$
\text{prox}_{\tau g}(\tilde{\bm{p}}) = \text{prox}_{\frac{\tau}{1+\tau\mu_g} \delta_{ \left\{\|\cdot\|_{2,\infty}\leq \lambda\right\} }}\left(\frac{\tilde{\bm{p}}}{1+\tau\mu_g}\right)= \Pi_{\left\{\|\cdot\|_{2,\infty}\leq \lambda\right\} }\left(\frac{\tilde{\bm{p}}}{1+\tau\mu_g}\right).
$$

Note that the same example has also been considered for similar verifications in \cite[Section 4.2]{Florea2016}: our results are in fact in good agreement with the ones reported therein.

\subsection*{Parameters} 
\begin{figure}[!h]
\begin{subfigure}[b]{0.31\textwidth}
\begin{center}
\includegraphics[height=3.5cm]{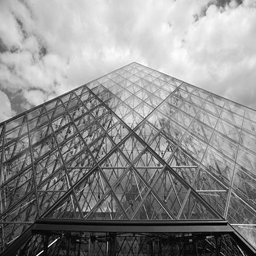}\vspace{0.1cm}
\end{center}
\caption{Original image}
\label{fig:original}
\end{subfigure}
\begin{subfigure}[b]{0.31\textwidth}
\begin{center}
\includegraphics[height=3.5cm]{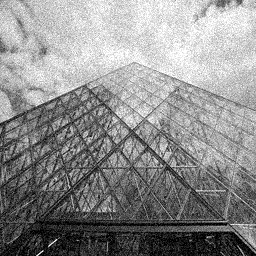}\vspace{0.1cm}
\end{center}
\caption{Noisy version}
\label{fig:noisy}
\end{subfigure}
\begin{subfigure}[b]{0.31\textwidth}
\begin{center}
\includegraphics[height=3.5cm]{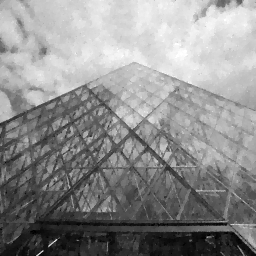}\vspace{0.1cm}
\end{center}
\caption{Denoised version}
\label{fig:denoised}
\end{subfigure}
\caption{Original, noisy and TV-Huber denoised images used. Noise is Gaussian distributed with zero mean and variance $\sigma^2=0.005$. The regularisation parameter is $\lambda=0.1$ and the Huber parameter is $\varepsilon=0.01$, so that $\mu=0.1$.}
\label{fig:denoising}
\end{figure}

In the following experiments we consider an image $u^0\in\R^{m\times n}$ with $m=n=256$ corrupted by Gaussian noise with zero mean and $\sigma^2=0.005$, see Figure \ref{fig:original}-\ref{fig:noisy}. We set the Huber parameter $\varepsilon=0.01$ and the regularisation parameter $\lambda=0.1$, so that $\mu_g=\mu=0.1$. In our comparisons we use the GFISTA algorithms \ref{alg:GFISTA} and \ref{alg:GFISTA:backtr} with and without backtracking using the prior knowledge of $L_f$ given by the estimate $L_f=8$ and an initial $L_0$, respectively. To ensure monotone decay we use the modified version described in Section \eqref{sec:monotone}, i.e. we use the modified update rules \eqref{monotone:update}-\eqref{monotone:update2}. For comparison, we report numerical results where the backtracking strategy is used `classically', i.e. it allows only for increasing of the Lipschitz constant estimate $L_k$ and used `adaptively' i.e. it allows for both its increasing and decreasing  along the iterations. The backtracking factor $\rho$ is set $\rho=0.9$. The initial value $t_1$ is set $t_1=1$.  The algorithm is initialised by the gradient of the noisy image $u^0$, i.e. $\bm{p_0}=Du^0$.

\begin{figure}[!h]
\begin{subfigure}[b]{0.49\textwidth}
\begin{center}
\includegraphics[width=6.2cm]{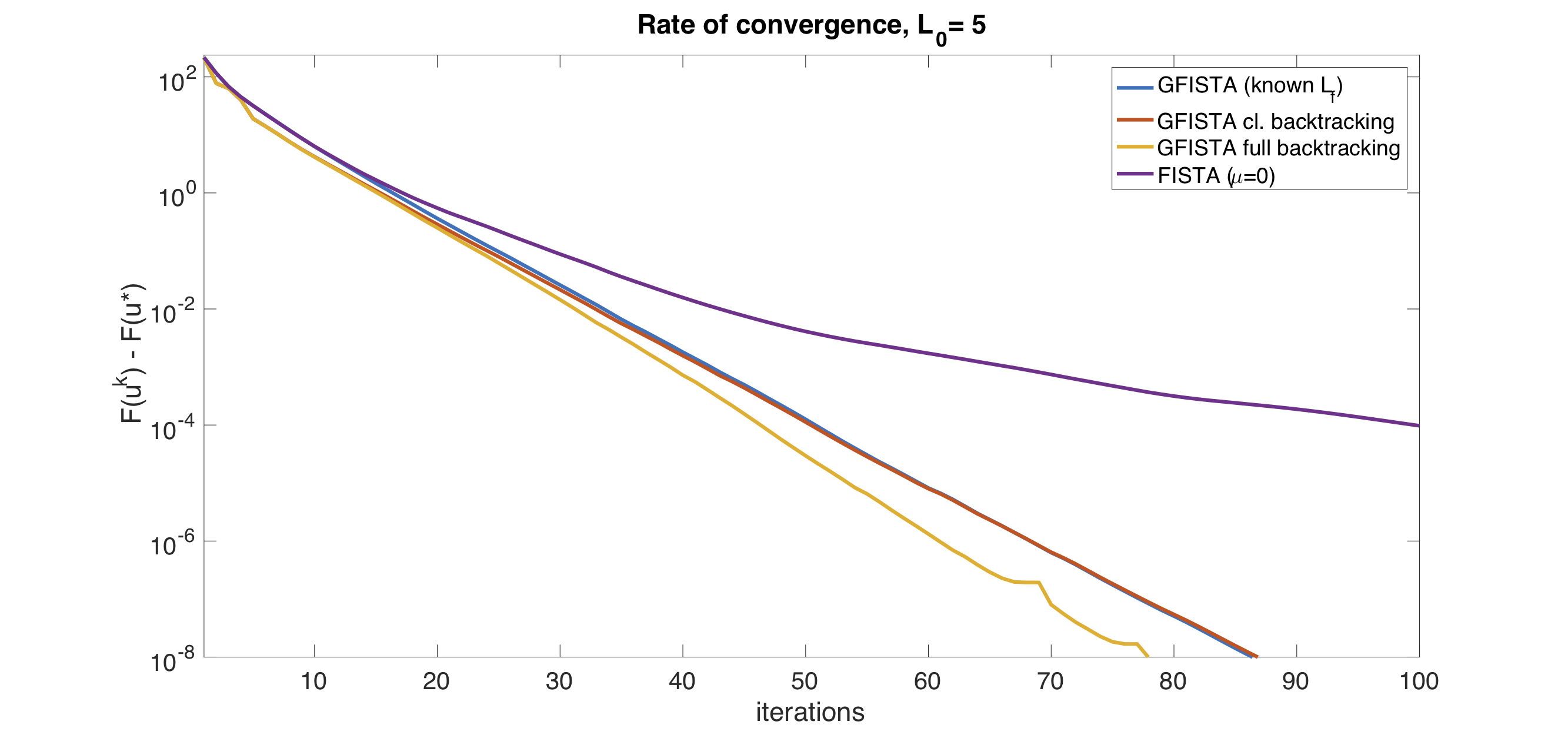}\vspace{0.1cm}
\end{center}
\caption{Convergence rates.}
\label{fig:rates:huber1:rates}
\end{subfigure}
\begin{subfigure}[b]{0.49\textwidth}
\begin{center}
\includegraphics[width=6.2cm]{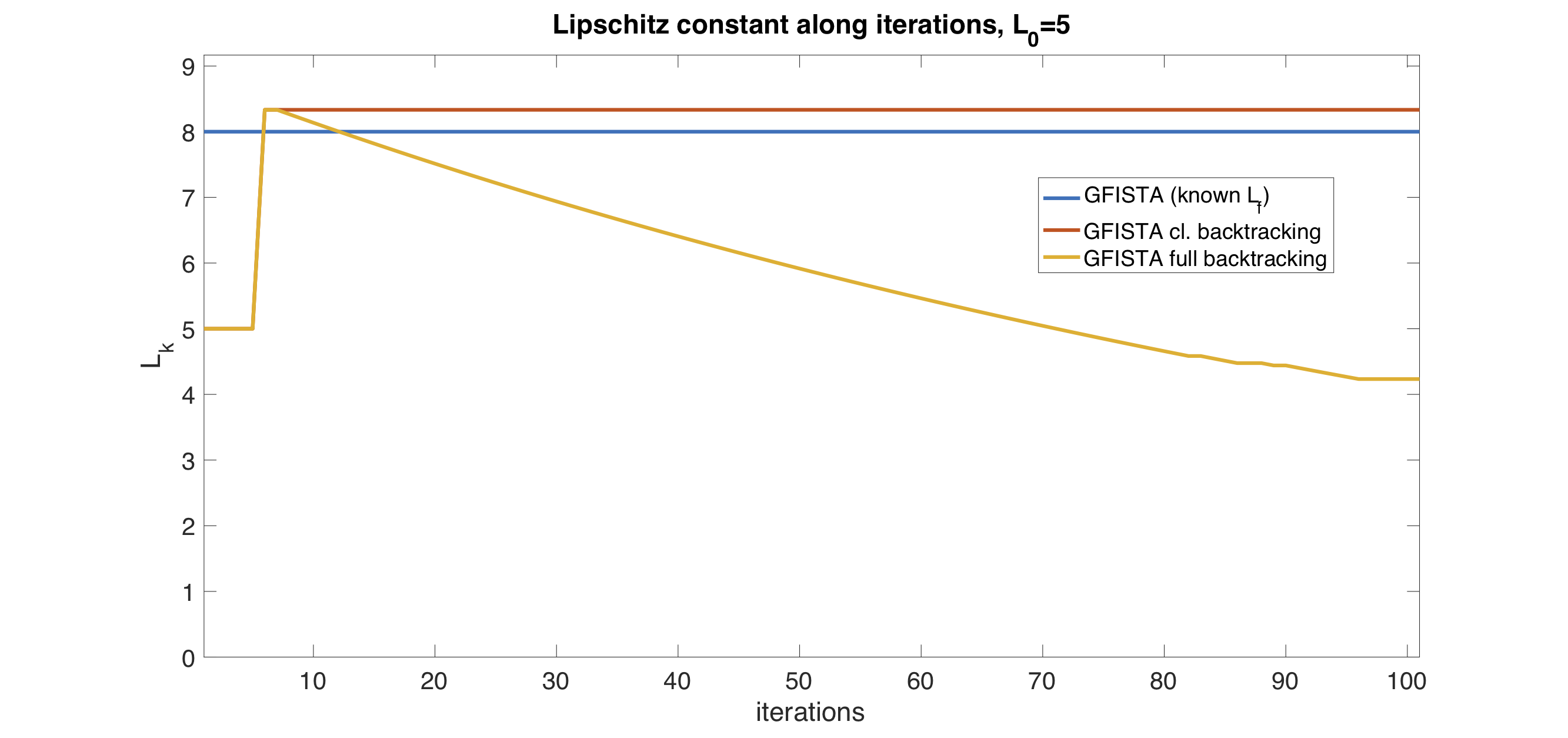}\vspace{0.1cm}
\end{center}
\caption{Lipschitz constant estimate.}
\label{fig:rates:huber1:LCE}
\end{subfigure}
\caption{Convergence rates and backtracking of the Lipschitz constant of $\nabla f$ starting from the underestimating initial value $L_0=5$.}
\label{fig:rates:huber1}
\end{figure}

\begin{figure}[!h]
\begin{subfigure}[b]{0.49\textwidth}
\begin{center}
\includegraphics[width=6.2cm]{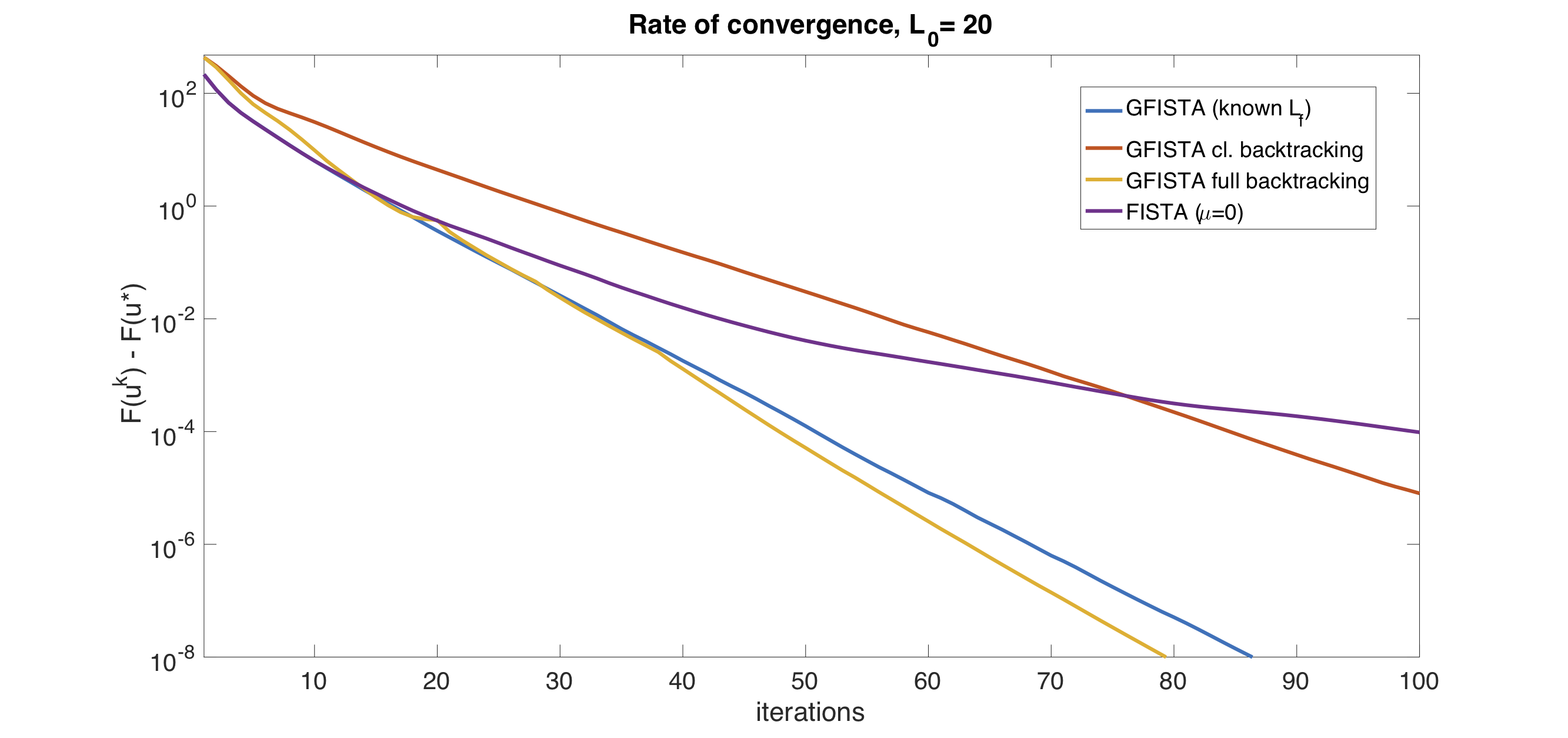}\vspace{0.1cm}
\end{center}
\caption{Convergence rates.}
\label{fig:rates:huber2:rates}
\end{subfigure}
\begin{subfigure}[b]{0.49\textwidth}
\begin{center}
\includegraphics[width=6.2cm]{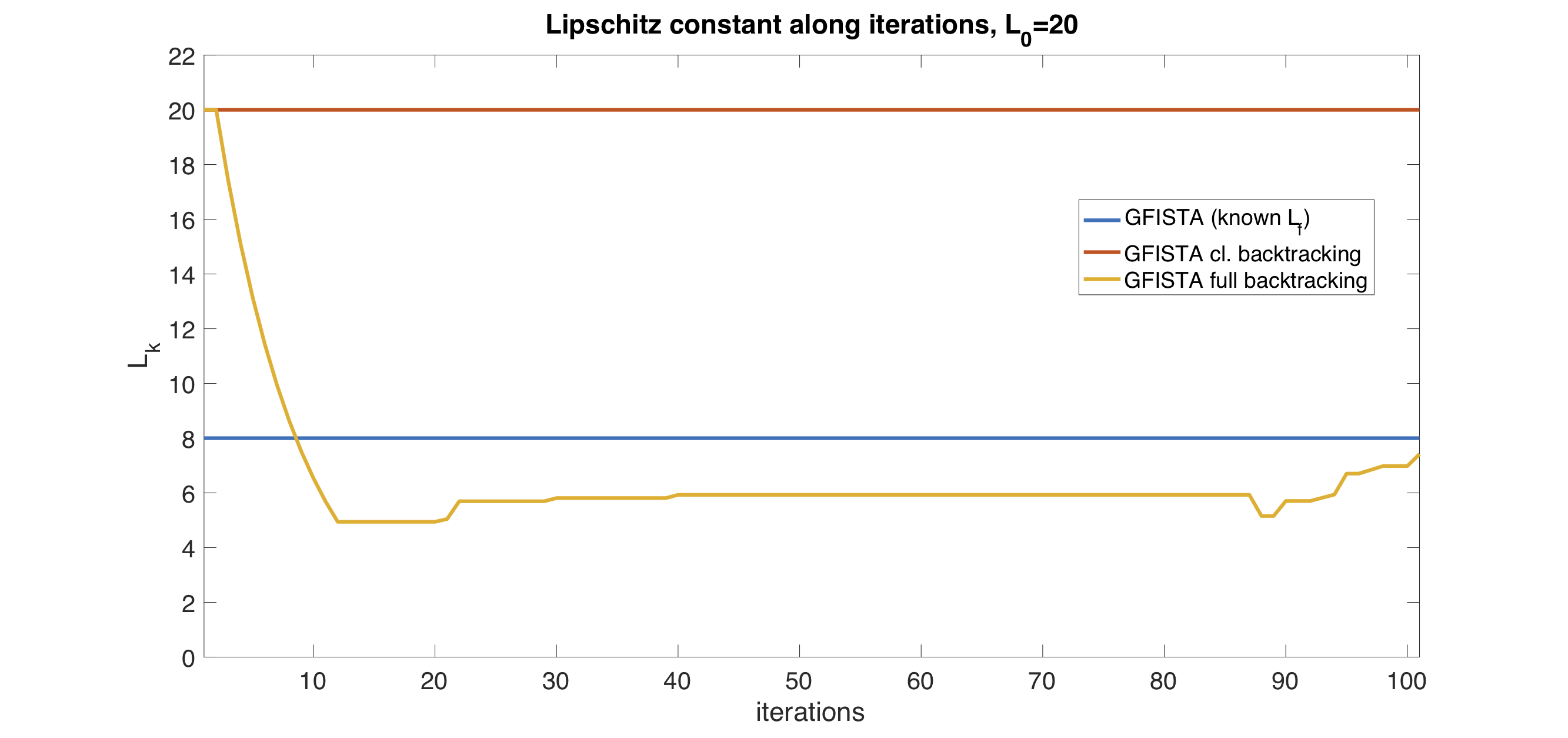}\vspace{0.1cm}
\end{center}
\caption{Lipschitz constant estimate.}
\label{fig:rates:huber2:LCE}
\end{subfigure}
\caption{Convergence rates and backtracking of the Lipschitz constant of $\nabla f$ starting from the overestimating initial value $L_0=20$.}
\label{fig:rates:huber2}
\end{figure}

To compute an approximation of the optimal solution $u^*$, we let the plain GFISTA algorithm run beforehand for $5000$ iterations and store the result for comparison, see Figure \ref{fig:denoised}. We then compute the results running the algorithms \ref{alg:GFISTA} and \ref{alg:GFISTA:backtr} for \texttt{iter}$=100$ iterations. We report the results computed for two different choices of $L_0$ which underestimate and overestimate the actual value of $L_f$, respectively, see Figure \ref{fig:rates:huber1} and \ref{fig:rates:huber2}.  For comparison, we further report the  $O(1/k^2)$ convergence rate of standard FISTA with no strongly convex parameter ($\mu=0$) encoded.

\subsection{Strongly convex TV Poisson denoising}  \label{sec:poisson:den}

In this second example we consider a different denoising model for images corrupted by Poisson noise, which is commonly observed in microscopy and astronomy imaging applications. Standard Poisson denoising models using Total Variation regularisation are typically combined with a convex,  non-differentiable Kullback-Leibler data fitting term, which can be consistently derived from the Bayesian formulation of the problem via MAP estimation (see, e.g., \cite{Saw11}). Here, we follow \cite{ChambolleEhrhardt2017} and consider a differentiable version of the Kullback-Leibler data term which, for a given positive noisy image $u^0\in\R^{m\times n}$ corrupted by  Poisson noise reads:
\begin{align}
f(u) & =\tilde{KL}(u_0,u)  \label{KL:diff} \\
&:=\sum_{i=1}^m\sum_{j=1}^n \begin{cases}
u_{i,j} + b_{i,j} - u^0_{i,j} + u^0_{i,j} \log\left( \frac{u^0_{i,j}}{u_{i,j}+b_{i,j}}\right)\qquad &\text{if }u_{i,j}\geq 0,\\
\frac{u^0_{i,j}}{2b_{i,j}^2}u_{i,j}^2 + \left(1-\frac{u^0_{i,j}}{b_{i,j}}\right)u_{i,j} + b_{i,j} - u^0_{i,j} + u^0_{i,j} \log\left( \frac{u^0_{i,j}}{b_{i,j}}\right)\qquad &\text{otherwise}, \notag
\end{cases}
\end{align}
where $b\in\R^{m\times n}$ stands for the background image which can be typically estimated from the data at hand. It is easy to verify the Lipschitz constant $\nabla \tilde{KL}(u_0,u)$  can be very roughly estimated as
 \begin{equation}  \label{est:Lipschitz}
 L_f= \max_{i,j} \frac{u^0_{i,j}}{b_{i,j}^2},
 \end{equation}
  which it is well-defined, positive and finite as long as $u^0$ and $b$ are positive. As a regularisation term, we will consider the following $\varepsilon$-strongly convex variant of isotropic TV in \eqref{discrete:TV}:
\begin{equation} \label{TV:L2:reg}
g(u)=\lambda\| Du\|_{2,1}+ \frac{\varepsilon}{2}\|u\|_2^2,
\end{equation}
where $\lambda>0$ stands again for the regularisation parameter. Differently from the Huber-TV ROF example, we aim here to apply the GFISTA algorithm \ref{alg:GFISTA:backtr} to solve composite problem:
\begin{equation}  \label{TV:strongly:poisson}
\min_u~ \lambda\| Du\|_{2,1}+ \frac{\varepsilon}{2}\|u\|_2^2 + \tilde{KL}(u_0,u)
\end{equation}
in primal form.

The gradient of the $KL$ term \eqref{KL:diff} can be easily computed and  the proximal map of $g$ in \eqref{TV:L2:reg} can be computed using the proximal map of the TV functional due to a general property reported in Lemma \ref{lemma:proximal} in the appendix, so that, recalling the definition \eqref{prox:metric}, for any $z$ there holds:
\begin{equation}  \label{prox:map}
\text{prox}_{\tau g}(z) = \text{prox}^{\frac{\lambda\tau}{1+\varepsilon\tau}}_{\|\cdot\|_{2,1}}\left(\frac{z}{1+\varepsilon\tau}\right).
\end{equation}
Thus, for any $\tau>0$, computing the right hand side of the equality above corresponds simply to solve the classical ROF problem with regularisation parameter $\sigma:=\frac{\lambda\tau}{1+\tau\varepsilon}$. We do that using standard FISTA as an iterative inner solver.

\subsection*{Parameters} 
\begin{figure}[!h]
\begin{subfigure}[b]{0.31\textwidth}
\begin{center}
\includegraphics[height=3.5cm]{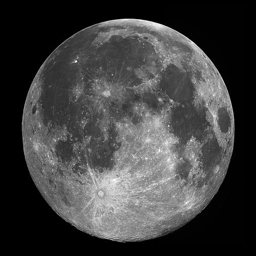}\vspace{0.1cm}
\end{center}
\caption{Original}
\label{fig:original_pois}
\end{subfigure}
\begin{subfigure}[b]{0.31\textwidth}
\begin{center}
\includegraphics[height=3.5cm]{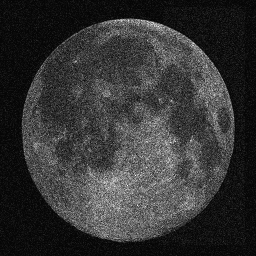}\vspace{0.1cm}
\end{center}
\caption{Noisy version}
\label{fig:noisy_pois}
\end{subfigure}
\begin{subfigure}[b]{0.31\textwidth}
\begin{center}
\includegraphics[height=3.5cm]{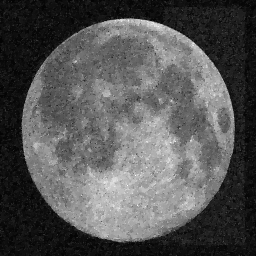}\vspace{0.1cm}
\end{center}
\caption{Denoised version}
\label{fig:denoised_pois}
\end{subfigure}
\caption{Original, noisy and restored image computed using the strongly convex TV-Poisson denoising model \eqref{TV:strongly:poisson}. The regularisation parameter is $\lambda=0.2$ and the strong convexity parameter is $\mu=\varepsilon=0.15$. }
\label{fig:denoising_pois}
\end{figure}

We consider an image $u^0\in\R^{m\times n}$ with $m=n=256$ corrupted artificially by Poisson noise, see Figure \ref{fig:original_pois}-\ref{fig:noisy_pois}. For simplicity, we consider a constant background with $b_{i,j}=1$ for all $i,j$. We set the strong convexity parameter $\varepsilon=0.15$ and the regularisation parameter $\lambda=0.1$. Clearly $\mu=\mu_g=\varepsilon$. In order to compute the proximal map \eqref{prox:map} we use $10$ iterations of standard FISTA. In the following example the Lipschitz constant of the gradient of the $\tilde{KL}$ term can be estimated via \eqref{est:Lipschitz} as $L_f=45$. We report in the following the results computed using the monotone variant of GFISTA algorithm \ref{alg:GFISTA} without backtracking and with classical and full backtracking (Algorithm \ref{alg:GFISTA:backtr} with monotone updates \eqref{monotone:update}-\eqref{monotone:update2}), for which the factor $\rho=0.8$ is chosen. The initial value $t_1$ is set $t_1=1$. The algorithm is initialised using the given noisy image $u^0$.

\begin{figure}[!h]
\begin{subfigure}[b]{0.49\textwidth}
\begin{center}
\includegraphics[width=6.2cm]{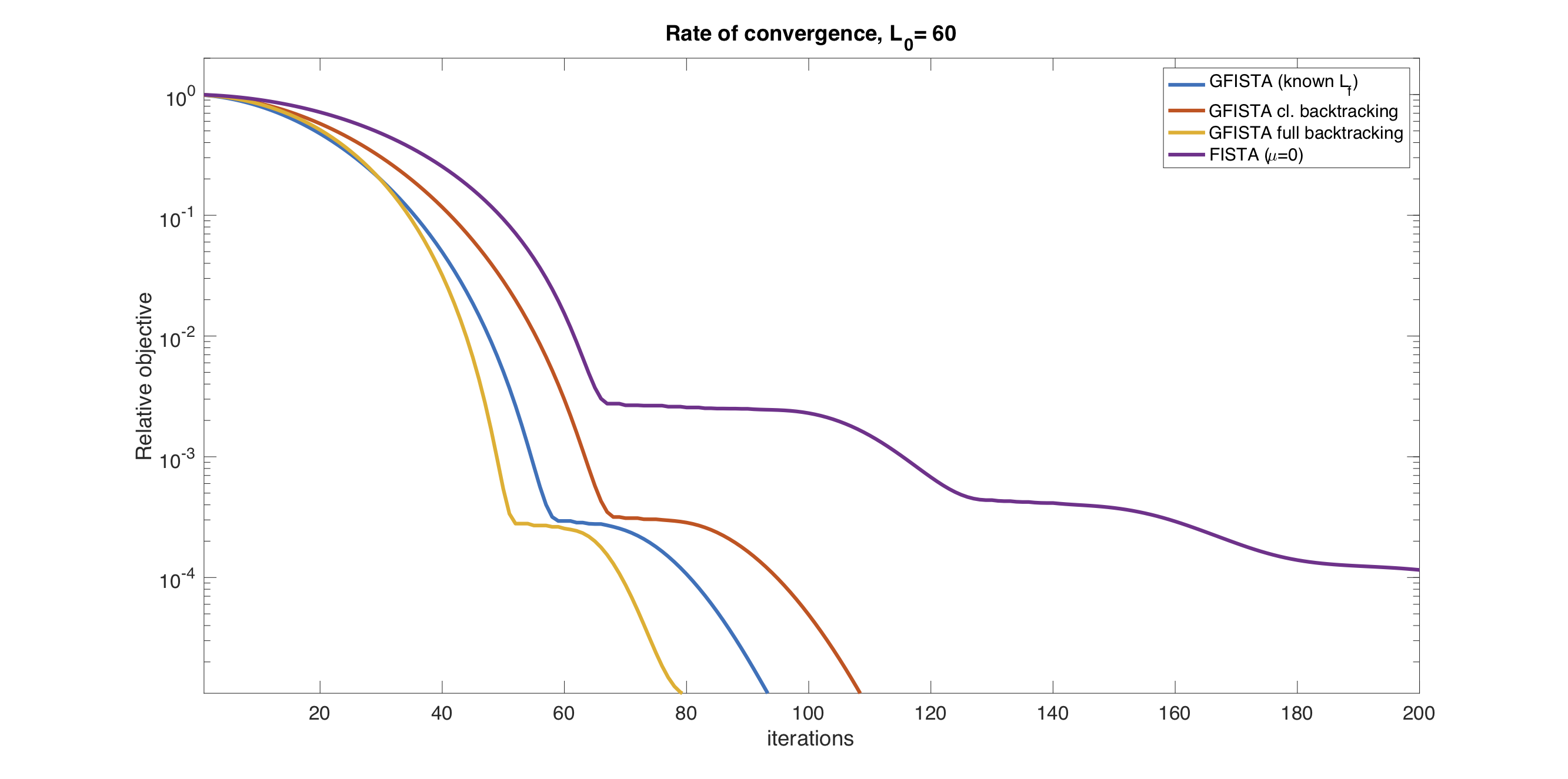}\vspace{0.1cm}
\end{center}
\caption{Convergence rates.}
\label{fig:rates:pois:rates}
\end{subfigure}
\begin{subfigure}[b]{0.49\textwidth}
\begin{center}
\includegraphics[width=6.2cm]{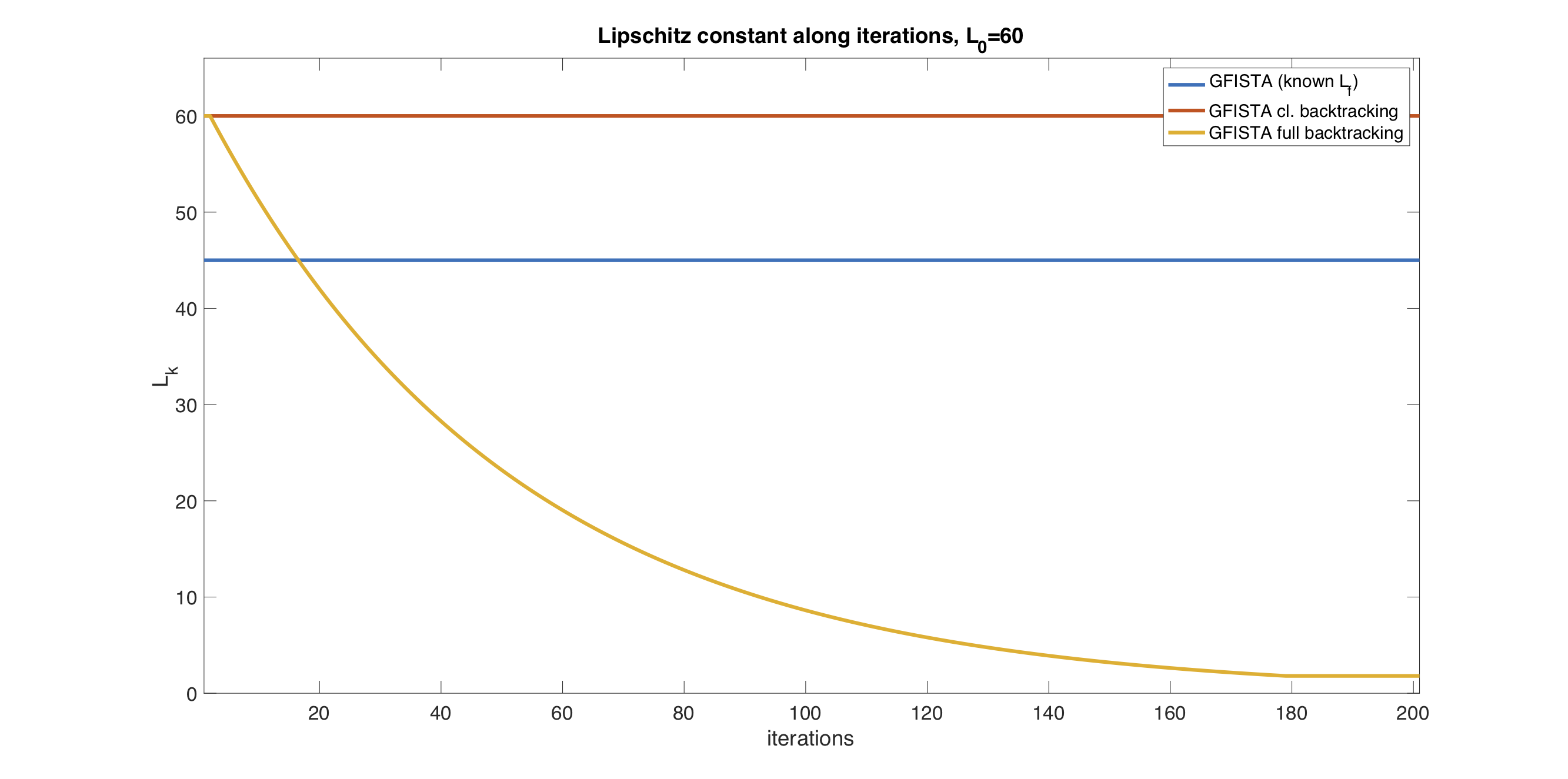}\vspace{0.1cm}
\end{center}
\caption{Lipschitz constant estimate.}
\label{fig:rates:pois:LCE}
\end{subfigure}
\caption{Convergence rates and backtracking of the Lipschitz constant of $\nabla f$ in \eqref{KL:diff} starting from the overestimating initial value $L_0=60$. Rates are shown in terms of the relative objective functional: $\frac{F(u^k)-F(u^*)}{F(u^0)-F(u^*)}$.}
\label{fig:rates:pois}
\end{figure}

\begin{figure}[!h]
\begin{center}
\includegraphics[height=4cm]{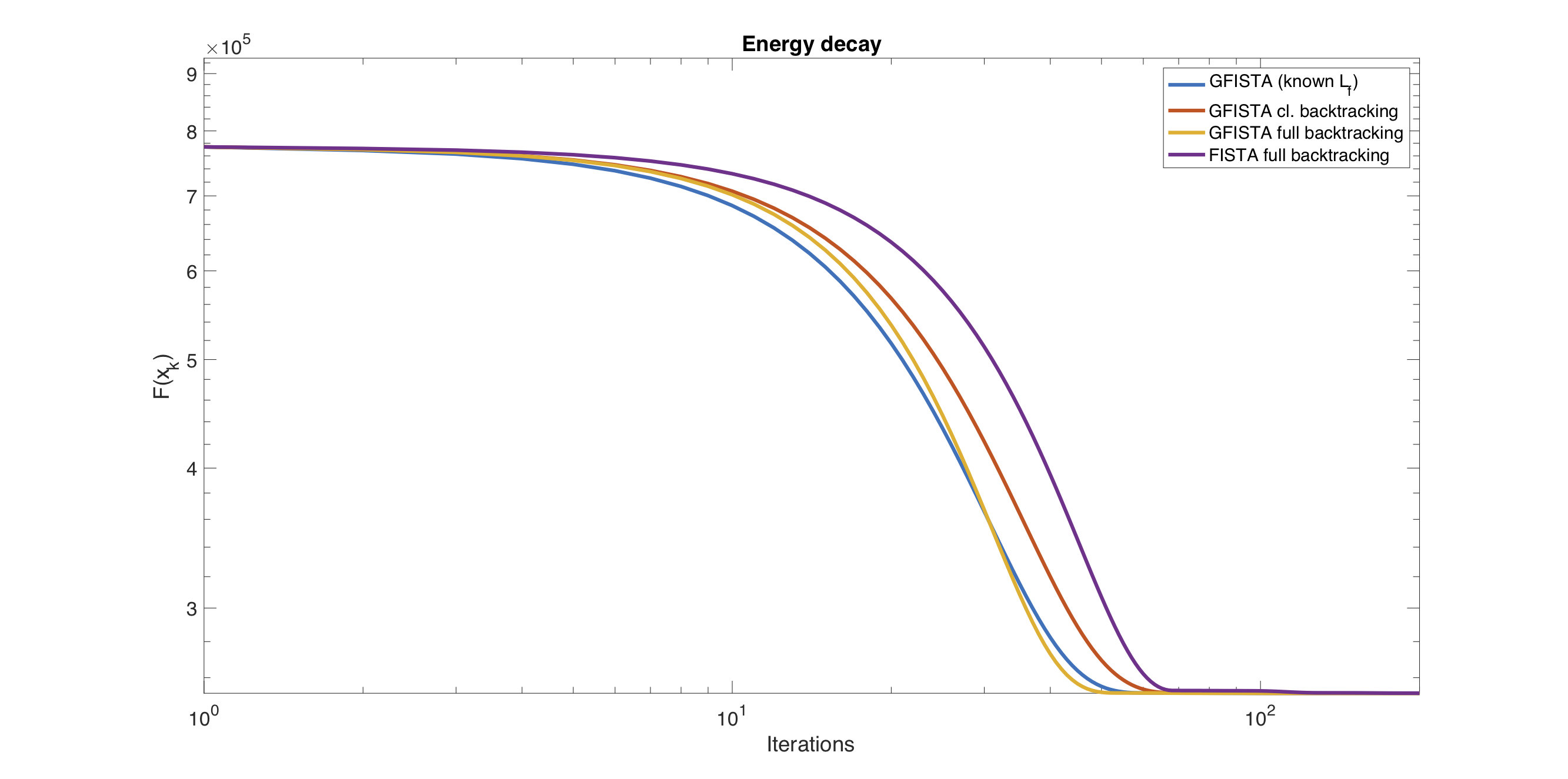}\vspace{0.1cm}
\end{center}
\caption{Monotone decay along the GFISTA iterates (with and without backtracking) after the monotone modification \eqref{monotone:update}-\eqref{monotone:update2}.}
\label{fig:energy_pois}
\end{figure}

An approximation of the solution $u^*$ is computed beforehand by letting the plain FISTA algorithm run for $5000$ iterations and then stored for comparison, see Figure \ref{fig:denoised_pois}. Results are then computed letting the monotone version of the GFISTA algorithms run for $\texttt{iter}=200$ iterations. In Figure \ref{fig:rates:pois} we report the results computed for a value of $L_0$ overestimating the actual one given by $L_f$ and in comparison with standard FISTA with no strongly convex modification. Once again we can observe that by incorporating the strongly convex modification of GFISTA linear convergence is achieved, in comparison with slower convergence of standard FISTA. Furthermore, the local estimate of the Lipschitz constant provided by the full backtracking strategy  decreases along the iterations, thus allowing for larger gradient steps and convergence in fewer iterations. In Figure \ref{fig:energy_pois}, we plot the monotone decay of the energy along the GFISTA iterates (with and without backtracking) after the monotone modification described in Section \eqref{sec:monotone}.

\subsection{Restarting strategies applied to the elastic net} \label{sec:restarting}
In this final example we test the performance of the GFISTA algorithm with backtracking \ref{alg:GFISTA:backtr}  in the case when a prior estimate of the strong convexity parameters $\mu_f$ and/or $\mu_g$ is either misspecified or not available. As a test problem we consider the Elastic Net regularisation model, which, for a given matrix $A\in\R^{m\times m}$, data $y\in\R^{m}$ and positive parameters $\lambda_1$ and $\lambda_2$ reads:
\begin{equation}  \label{eq:Elastic_Net}
\min_u~\left\{  F(u):=\frac{1}{2}\| Au-y\|_2^2 + \lambda_1\|u\|_1 + \frac{\lambda_2}{2}\|u\|_2^2  \right\},
\end{equation}
The Elastic Net is commonly used in the study of logistic regression models as a regularised version of the LASSO estimator by means of a ridge-type quadratic term and it is employed for several parameter identification \cite{Zou2005} and support vector machine problems \cite{Wang2006}.
In order to apply the the GFISTA algorithm \ref{alg:GFISTA:backtr}, we split the functional $F$ above into the sum:
\begin{equation}   \label{eq:f_g_elastic_net}
f(u):= \frac{1}{2}\| Au-y\|_2^2 + \frac{\lambda_2}{2}\|u\|_2^2 , \qquad g(u):= \lambda_1\|u\|_1.
\end{equation}
Under this choice, we note that $f$ is differentiable with Lipschitz-continuous gradient given by
$$
\nabla f(u)= A^*(Au-y) + \lambda_2 u
$$
 whose Lipschitz constant can be calculated as $L_f = \lambda_{max}(A^*A + \lambda_2\bm{Id})$, where by $\lambda_{max}(M)$ we denote the largest eigenvalue of the matrix $M$. Note that in case of large-size problems ($m\gg 1$), such computation of $L_f$ may render prohibitively expensive. The non-smooth function $g$ is convex and for $\tau>0$ its proximal map can be calculated component-wise by the soft-thresholding operator as:
$$
\Big(\mathrm{prox}_{\tau g}(z)\Big)_i= \mathrm{sign}(z_i)\max\left(|z_i|-\tau\lambda_1,0 \right),\quad i=1,\ldots,m.
$$
Finally, note that $f$ is $\lambda_2$-strongly convex, so that $\mu=\mu_f=\lambda_2$.
\subsection*{Parameters} In the following experiments we solve the problem \eqref{eq:Elastic_Net} in correspondence of a normalised randomly generated operator $A\in\R^{3600\times 3600}$ and for parameters $\lambda_1, \lambda_2$ set as $\lambda_1=0.01$ and $\lambda_2=1e^{-5}$, so that $\mu=\mu_f=\lambda_2$. The Lipschitz constant $L_f$ of $\nabla f$ can be estimated in this example as $L_f=0.0657$. For the backtracking routine, we set the backtracking factor $\rho=0.95$. The GFISTA algorithm \ref{alg:GFISTA:backtr} is initialised by $t_1=1$, $L_0=1$ and $x_0=\bm{0}$.  The plain GFISTA algorithm \eqref{alg:GFISTA} without backtracking is run for $5000$ iterations and its solution $x^*$ is stored for comparisons. The following results are computed by running the algorithm for $\texttt{iter}=100$ iterations.

\begin{figure}[!h]
\begin{subfigure}[b]{0.49\textwidth}
\begin{center}
\includegraphics[width=6.2cm]{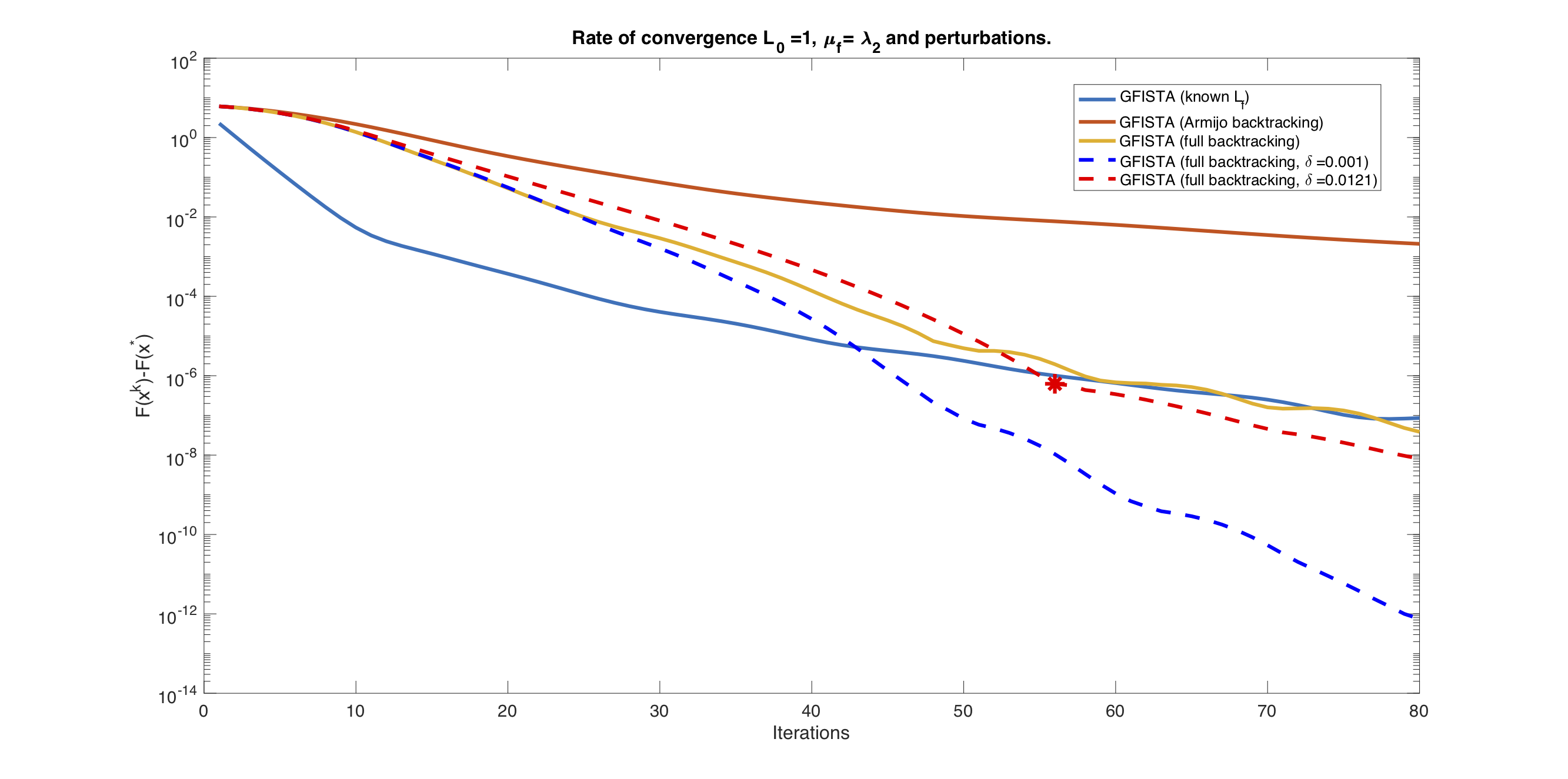}\vspace{0.1cm}
\end{center}
\caption{Convergence rates.}
\label{fig:rates:ENET:rates}
\end{subfigure}
\begin{subfigure}[b]{0.49\textwidth}
\begin{center}
\includegraphics[width=6.2cm]{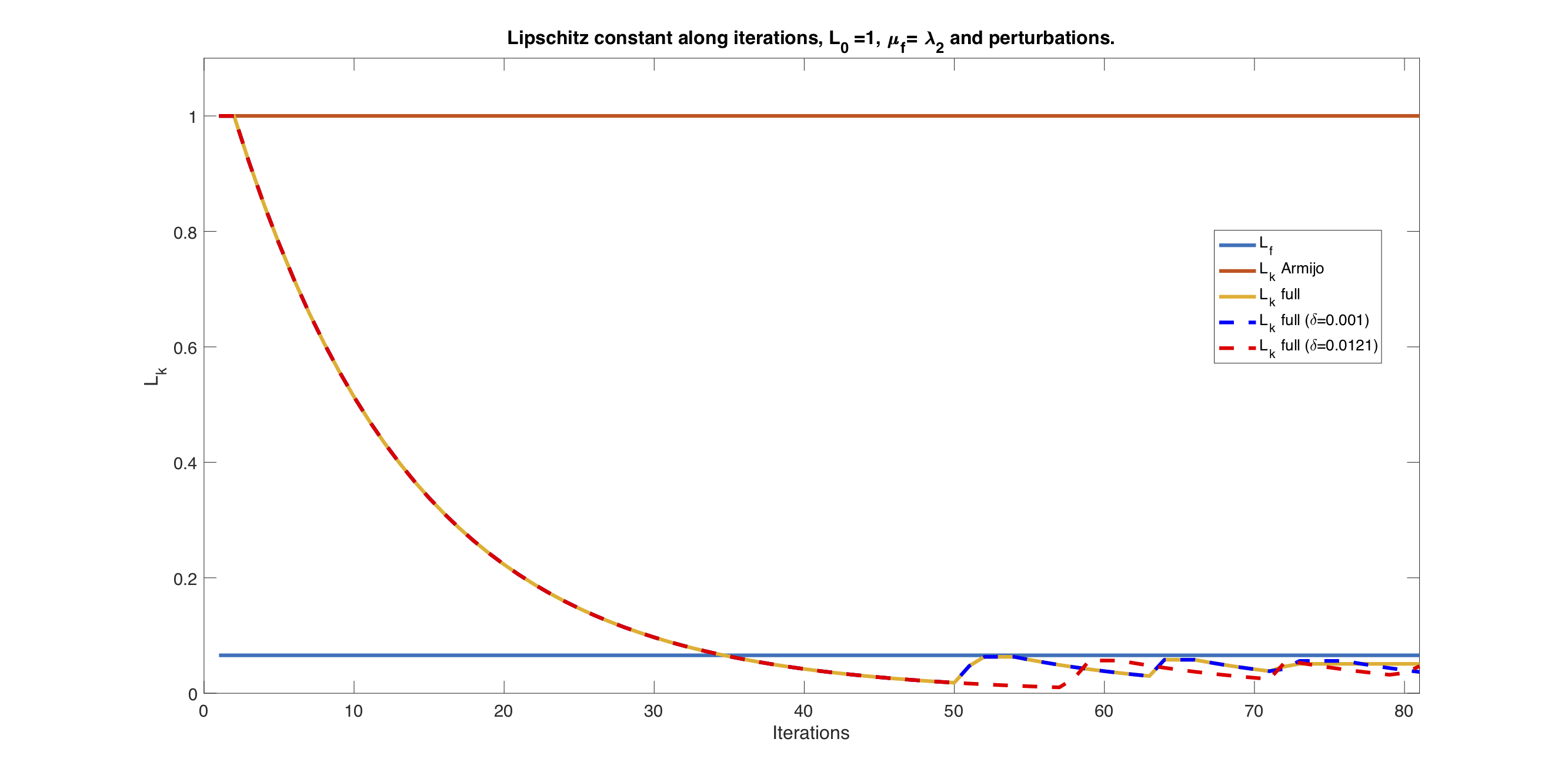}\vspace{0.1cm}
\end{center}
\caption{Lipschitz constant estimate.}
\label{fig:rates:ENET:LCE}
\end{subfigure}
\caption{Convergence rates and backtracking of the Lipschitz constant of $\nabla f$ in \eqref{eq:f_g_elastic_net} starting from the overestimating initial value $L_0=1$. In the convergence plot, full lines refer to the case when $\mu_f=\lambda_2$, while dashed lines refer to ``wrong" values of $\mu_f$ which is perturbed as $\mu_f=\lambda_2+ \delta$. In the red dashed line we also scatter the point corresponding to the iteration violating the condition $\mu_f< L_k$ along the iterations, which requires the reduction of $\mu_f$.}
\label{fig:rates:ENET}
\end{figure}

In the first test, we compare once again the performance of the GFISTA algorithm \ref{alg:GFISTA:backtr} when the prior estimate of $L_f$ is available and when it is not, using both standard Armijo-type backtracking and the adaptive one proposed in this work, see Figure \ref{fig:rates:ENET}. Compared to the examples considered above, note that in this case the strong convexity constant of the problem is encoded in the term $f$ defined in \eqref{eq:f_g_elastic_net}, which is accommodated by our strategy. Note, however, that it renders typically more efficient to encode strong convexity in the non-smooth component $g$ which is treated implicitly rather than in $f$ which is treated explicitly. This latter choice would require in fact more restrictive time-steps $\tau\leq 1/(L_f + \mu_f)$. 

\noindent In addition, we also report the results obtained when a ``wrong" value of $\mu_f$ is used. Given its quadratic behaviour, one may in fact suppose that in addition to the $\lambda_2$-strongly convexity, some further strong convexity could be hidden in the quadratic data fitting term. In the following, we then report the results obtained by applying the GFISTA algorithm \ref{alg:GFISTA:backtr} with full backtracking for a perturbed value of $\mu_f $ given by $\mu_f= \lambda_2 + \delta$, for a small perturbation $0< \delta\ll 1$.  Note that under such modification the natural condition $\mu_f< L_k$ may be violated along the iterations, thus preventing the algorithm from converging. Whenever this happens, we decrease the value $\mu_f$ of a factor $\rho$, redefine the term $q_k$ appearing in Algorithm \ref{alg:GFISTA:backtr} in correspondence of this new value and carry on with the algorithm. In this way convergence is always guaranteed and also large misspecifications of $\mu_f$ can be treated.

\noindent Provided such verification is performed along the iterations, these tests suggest that encoding further, hidden, strong convexity information in the model \eqref{eq:Elastic_Net} can improve the convergence rates of the GFISTA algorithm \ref{alg:GFISTA:backtr}.  

\begin{figure}[!h]
\begin{subfigure}[b]{0.49\textwidth}
\begin{center}
\includegraphics[width=6.2cm]{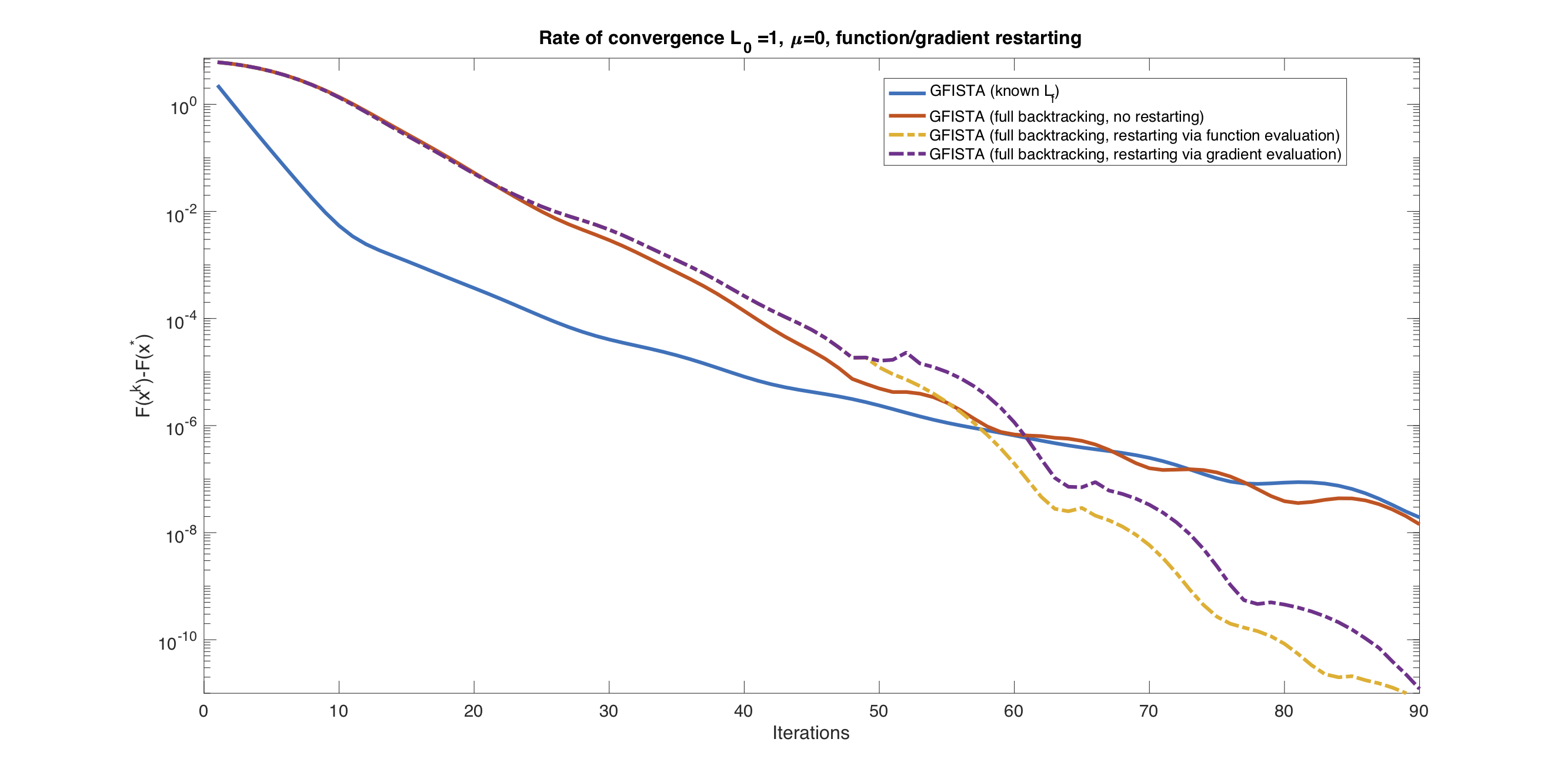}\vspace{0.1cm}
\end{center}
\caption{Convergence rates.}
\label{fig:rates:ENET1:rates}
\end{subfigure}
\begin{subfigure}[b]{0.49\textwidth}
\begin{center}
\includegraphics[width=6.2cm]{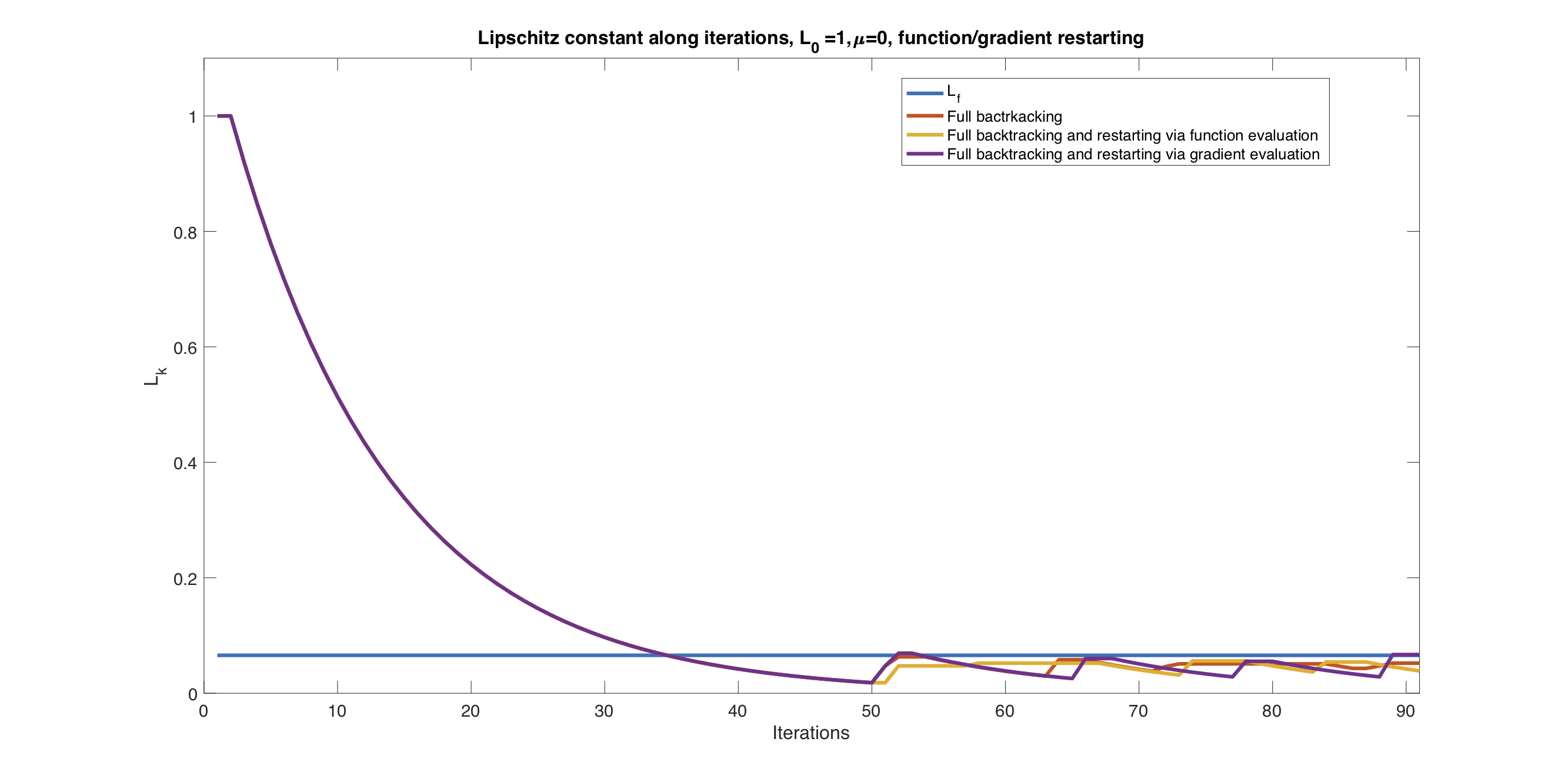}\vspace{0.1cm}
\end{center}
\caption{Lipschitz constant estimate.}
\label{fig:rates:ENET1:LCE}
\end{subfigure}
\caption{Convergence rates and backtracking of the Lipschitz constant of $\nabla f$ in \eqref{eq:f_g_elastic_net} with and without restarting based on function \eqref{eq:restarting_function} and gradient  \eqref{eq:restarting_gradient} criteria. Initial overestimating value $L_0=1$.}
\label{fig:rates:ENET1}
\end{figure}

Motivated by these considerations, we perform in the following a further numerical test where we assume that the values of the strong convexity parameters $\mu_f$ and $\mu_g$ (and, consequently, $\mu$) are unknown. In several applications, it is actually very hard to provide an explicit estimation of such parameters.  Moreover, as we have seen in the examples above, some hidden strong convexity can be still not detected explicitly only looking at the structure of the functions $f$ and $g$. An indirect way to estimate strong convexity consists in restarting the algorithm depending on a certain criterion, see, e.g., \cite{Nesterov2013}. In \cite{Candes2015} two heuristic restarting procedures based either on  the evaluation of the composite functional or of a (generalised) gradient are studied. These two restarting approaches have become very popular since then and, more recently, some others have been proposed, for instance in \cite{Lin2015} and \cite{Ferocq2016}. Here, we follow \cite{Candes2015} and apply the two function- and gradient-based restarting procedures to the GFISTA algorithm \ref{alg:GFISTA:backtr} with full backtracking to solve the Elastic Net problem above under the same choice of parameters as above.
As discussed in \cite[Section 5.2]{Candes2015} the two restarting criteria to consider for FISTA-type algorithms are the following:
\begin{itemize}
\item \textbf{Function adaptive restart}: restart the algorithm whenever 
\begin{equation} \label{eq:restarting_function}
F(u^{k+1})>F(u^{k}).
\end{equation}
\item \textbf{Gradient adaptive restart}: restart the algorithm whenever 
\begin{equation} \label{eq:restarting_gradient}
(y^k-u^{k+1})^T(u^{k+1}-u^k)>0.
\end{equation}
\end{itemize}
Compared to the function-based restarting scheme, the gradient adaptive restart is  observed to be more stable around $x^*$. Furthermore, there is no extra computational cost in applying such restarting to GFISTA  \ref{alg:GFISTA:backtr} since all the quantities appearing in \eqref{eq:restarting_gradient} have already been calculated during the backtracking phase. We remark that this second approach goes under the name of `gradient' restart since one can interpret for each $k\geq 0$ the FB step \eqref{rules:update:xk1} in Algorithm \ref{alg:GFISTA:backtr} as a \emph{generalised} gradient step in defined by
$$
x^{k+1} =  \text{prox}_{\tau_{k+1} g} (y^{k} - \tau_{k+1}\nabla f(y^{k}))  =: y^k - \tau_{k+1} G(y^k).
$$
The restarting condition \eqref{eq:restarting_gradient} would then actually read in this case $G(y^k)^T(u^{k+1}-u^k)>0$.
In Figure \ref{fig:rates:ENET1}, we report the convergence plots and the Lipschitz constant variations for the solution of the Elastic Net problem \eqref{eq:Elastic_Net} via the GFISTA algorithm \ref{alg:GFISTA:backtr} with full backtracking combined with the two restarting strategies above. We observe a faster linear convergence compared to the fully backtracked GFISTA algorithm which, heuristically, can therefore be adapted and efficiently employed also to strongly convex problem with no prior estimate on the strong convexity constant $\mu$. A rigorous proof of these convergence results is left for future research.

\section{Conclusions and outlook}  \label{sec:conclusions}

We study a fast backtracking strategy for the strongly convex variant of the FISTA algorithm proposed in \cite{ChambollePock2016} and based on a inequality condition expressed in terms of the Bregman distance, see Section \ref{FISTA-CP}. Using standard properties of strongly convex functions and upon multiplication by appropriate terms, we have derived in Section \ref{sec:backtr_GFISTA} the convergence estimate \eqref{decay5} whose decay factor \eqref{def:theta_k} has been then studied carefully to estimate the convergence speed of Algorithm \ref{alg:GFISTA:backtr}. Our analysis is essentially based on classical technical tools similar to the ones used in Nesterov in \cite{Nesterov2004} and on general properties of the extrapolation sequences defined. Our main result is reported in Theorem \ref{theo:convergence} where accelerated linear convergence rates are proved in term of average quantities depending on the estimated values along the iterations.  Our theoretical results are verified numerically in Section \ref{sec:numerics} on some exemplar problems.

The backtracking strategy proposed is fast and robust since it allows for adaptive adjustment of the gradient step size (i.e. the proximal map parameter) depending on the local `flatness' of the gradient of the component $f$ in the objective functional, i.e. on the local estimate $L_k$ of $L_f$. In other words, in flat regions (small $L_f$) larger step sizes are promoted, whereas where large variations of $\nabla f$ occur (large $L_f$), small steps are preferred for a more accurate descent. From an algorithmic point of view, extrapolation is performed using suitable parameters providing strict decay in the convergence inequality \eqref{decay5} and defined not only in terms of the step sizes, but also in terms of the strong convexity parameters of $f$ and $g$ and resulting in more refined convergence rate estimates. 
Finally, in terms of computational costs our approach has a lower per-iteration cost than the one studied by Nesterov in \cite{Nesterov2013} since it avoids the calculation of the gradient of the smooth component in the proximal step. Accelerated convergence rates are proved and defined in terms of average quantities depending on the estimates performed along the iterations.

Further research could address the rigorous analysis of the combined backtracking approach with the restarting procedures \'a la Cand\'es used in Section \ref{sec:restarting} for situations when  the strong convexity parameters $\mu_f$ and $\mu_g$ are unkonwn. In this work we have heuristically showed good performance only for the case of function- and gradient-based restarting procedures, but it would be of great interest also exploring more the recently proposed approaches by Fercoq and Qu \cite{Ferocq2016} where the restarting does not require any condition but combines appropriately past iterates of the algorithm in an appropriate way. A rigorous analysis of a combined backtracking-restarting procedure would be very interesting for the sake of designing an algorithm fully adaptive to local convexity and smoothness of its functions.

Finally, it would be interesting to test the robustness and the performance of our algorithm on other strongly convex, possibly large-scale problems coming from the fields of image and data analysis with various condition numbers.

\section*{Acknowledgements}
The authors would like to  thank the valuable comments of the anonymous referees which improved significantly the quality of the manuscript.

\appendix

\section{Some useful Lemmas}   \label{appendix} 

In this appendix we prove some general results which has been used in our work. We start with a general inequality used to derive the descent rule \eqref{descent:rule}. Its proof is a consequence of a trivial property of strongly convex functions.

\newtheorem{notation}{Lemma}[section]

\begin{notation}
If $h:\mathcal{X}\to\R\cup\left\{\infty\right\}$ is strongly convex with parameter $\mu_h>0$ and $\hat{x}\in\mathcal{X}$ is a minimiser of $h$, the following property holds:
\begin{equation}  \label{prop:strong:conv}
h(x)\geq h(\hat{x}) + \frac{\mu_h}{2}\|x-\hat{x}\|^2,
\end{equation}
for any $x\in\mathcal{X}$.
\end{notation}

\begin{proof}
By definition of $\mu_h$-strong convexity, for any $x,y\in\mathcal{X}$ there holds:
$$
h(x) \geq h(y) + \langle p, y-x\rangle + \frac{\mu_h}{2}\|x-y\|^2,
$$
where $p\in \partial h(y)$, the subdifferential of $h$ evaluated in $y$.
Taking $y=\hat{x}$, since $0\in \partial h(\hat{x})$, we get \eqref{prop:strong:conv}.
\end{proof}

An immediate consequence of this general property is the proof of the descent rule \eqref{descent:rule} used in Section \ref{FISTA-CP} as a starting point of our convergence estimates. We follow \cite{ChambollePock2016,Tseng2008}.

\newtheorem{lemma:descent}[notation]{Lemma}

\begin{lemma:descent}   \label{lemma:descent:proof}
Let $f:\mathcal{X}\to\R$ be a $\mu_f$-strongly convex function with Lipschitz gradient with constant $L_f$ and $g:\mathcal{X}\to\R\cup\left\{\infty\right\}$ be a l.s.c., $\mu_g$-strongly convex function. Then, defining for any $\bar{x}\in\mathcal{X}$ and any $0<\tau <1/L_f$ the forward-backward map: $T_\tau: \bar{x}\mapsto \text{prox}_{\tau g}\left(\bar{x}-\tau\nabla f(\bar{x})\right)=:\hat{x}$, the following inequality holds for the composite functional $F=f+g$:
\begin{equation} \label{descent:rule:app}
F(x)+(1-\tau\mu_f)\frac{\|x-\bar{x}\|^2}{2\tau} \geq F(\hat{x}) + (1+\tau\mu_g)\frac{\|x-\hat{x}\|^2}{2\tau}, \qquad \text{ for any } x\in \mathcal{X}.
\end{equation}
\end{lemma:descent}

\begin{proof}
By definition, $\hat{x}$ is the minimiser of the function $h:\mathcal{X}\to\R\cup\left\{\infty\right\}$ defined by:
$$
h:x\mapsto g(x) + f(\bar{x}) + \langle f(\bar{x},x-\bar{x}\rangle + \frac{\|x-\bar{x} \|^2}{2\tau}.
$$
The function $h$ is strongly convex with parameter $\mu_h:=(\tau\mu_g+1)/\tau$. Hence, for any $x\in \mathcal{X}$:
\begin{align}
& F(x)+(1-\tau\mu_f)\frac{\|x-\bar{x}\|^2}{2\tau} \geq g(x) + f(\bar{x}) + \langle \nabla f(\bar{x}),x-\bar{x}\rangle + \frac{\|x-\bar{x} \|^2}{2\tau} \notag\\
& \geq g(\hat{x}) + f(\bar{x}) + \langle \nabla f(\bar{x}),\hat{x}-\bar{x}\rangle + \frac{\|\hat{x}-\bar{x} \|^2}{2\tau} + (1+\tau\mu_g)\frac{\|x-\hat{x}\|^2}{2\tau}\notag\\
& \geq g(\hat{x}) + f(\hat{x}) + \frac{1-\tau L_f}{2\tau}\| \hat{x}-\bar{x}\|^2 + (1+\tau\mu_g)\frac{\|x-\hat{x}\|^2}{2\tau}, \notag \\
& = F(\hat{x}) + \frac{1-\tau L_f}{2\tau}\| \hat{x}-\bar{x}\|^2 + (1+\tau\mu_g)\frac{\|x-\hat{x}\|^2}{2\tau}, \label{RHS:lemma}
\end{align}
where the first inequality holds by strong convexity of $f$, the second one is a simple application of Lemma \ref{prop:strong:conv} and the last one follows from the Lipschitz continuity of $\nabla f$. Since $\tau L_f<1$ by assumption, we can neglect the third term in \eqref{RHS:lemma} and get \eqref{descent:rule:app}.
\end{proof}

We finally report a general properties of proximal mappings which we used in our numerical experiments in Section \ref{sec:numerics}. For a general convex function $h$ it essentially allows a straightforward calculation of the proximal map of the composite $\varepsilon$-strongly convex function $g:=\alpha h + \frac{\varepsilon}{2}\|\cdot\|_2^2$ in terms of the proximal map of $h$ itself. We recall the notation \eqref{prox:metric}.

\newtheorem{lemma:proximal:map}[notation]{Lemma}

\begin{lemma:proximal:map}  \label{lemma:proximal}
Let $h:\mathcal{X}\to\R\cup\left\{+\infty\right\}$ a convex, proper and l.s.c. function. For $\alpha, \varepsilon>0$ let $g$ be defined as:
$$
g(x):=\alpha h(x) + \frac{\varepsilon}{2}\|x\|^2,\qquad x\in \mathcal{X}.
$$
Then, there holds:
\begin{equation*} 
\mathrm{prox}_{\tau g}(z) = \mathrm{prox}_{h}^{\frac{\alpha\tau}{1+\varepsilon\tau}}\left(\frac{z}{1+\varepsilon\tau}\right),\quad \text{for any }\tau>0\text{ and }z\in\mathcal{X}.
\end{equation*}
\end{lemma:proximal:map}

\begin{proof}
Let $\tau>0$ and $z\in\mathcal{X}$. We have the following chain of equalities:
\begin{align}
& \mathrm{prox}_{\tau g}(z)  =\mathrm{prox}^\tau_{g}(z) = \argmin_{y\in\mathcal{X}} ~g(y) + \frac{1}{2\tau}\|y-z\|^2  \notag\\
& =\argmin_{y\in\mathcal{X}}~ h(y) + \frac{1+\tau\varepsilon}{2\alpha\tau}\|y\|^2 + \frac{1}{2\alpha\tau}\|z\|^2 - \frac{1}{\alpha\tau}\langle y, z\rangle \notag\\
& =\argmin_{y\in\mathcal{X}}~ h(y) + \frac{1}{2\frac{\alpha\tau}{1+\tau\varepsilon}}\|y\|^2 + \Big(\frac{1}{2\alpha\tau(1+\varepsilon\tau)}- \frac{\varepsilon}{2\alpha(1+\varepsilon\tau)}\Big)\|z\|^2 - \frac{1}{\alpha\tau}\langle y, z\rangle \notag \\
& =\argmin_{y\in\mathcal{X}}~ h(y) + \frac{1}{2\frac{\alpha\tau}{1+\tau\varepsilon}}\|y\|^2 +  \frac{1}{2\frac{\alpha\tau}{1+\tau\varepsilon}} \|\frac{z}{1+\varepsilon\tau}\|^2 - \frac{1+\varepsilon\tau}{\alpha\tau}\langle y, \frac{z}{1+\varepsilon\tau}\rangle \notag \\
& = \argmin_{y\in\mathcal{X}}~ h(y) +  \frac{1}{2\frac{\alpha\tau}{1+\tau\varepsilon}} \| y - \frac{z}{1+\varepsilon\tau}\|^2= \mathrm{prox}^{\frac{\alpha\tau}{1+\varepsilon\tau}}_{h}(\frac{z}{1+\varepsilon\tau}). \notag
\end{align}
\end{proof}

\bibliographystyle{amsplain}
\bibliography{FISTAbib}
\end{document}